\documentclass[a4paper,11pt]{amsart}
\usepackage[utf8]{inputenc}
\usepackage[foot]{amsaddr}

\usepackage[english]{babel}
\usepackage{amsmath,amssymb,verbatim,mathrsfs,latexsym,paralist}
\usepackage{graphicx}
\usepackage{epstopdf}
\usepackage{indentfirst}
\usepackage{multirow}
\usepackage{amsthm}
\allowdisplaybreaks[4]
\usepackage[pdftex]{hyperref}

\newcommand{\Z}{\mathbb{Z}}
\newcommand{\C}{\mathbb{C}}

\newcommand{\mf}{\mathfrak}
\newcommand{\g}{\mf{g}}
\newcommand{\h}{\mf{h}}

\newcommand{\Hom}{\mathrm{Hom}}

\numberwithin{equation}{section}
\newtheorem{theorem}{Theorem}[section]
\newtheorem*{theorem*}{Theorem}
\newtheorem{thm}{Theorem}[section]
\newtheorem{proposition}[theorem]{Proposition}
\newtheorem{conj}[theorem]{Conjecture}
\newtheorem{lemma}[theorem]{Lemma}
\newtheorem{corollary}[theorem]{Corollary}
 
\newtheorem{defn}[theorem]{Definition}            
\theoremstyle{remark}
 
\theoremstyle{remark}
\newtheorem{rmk}[theorem]{Remark}
\newtheorem{example}[theorem]{Example}   

\title[On the Spectrum of Exterior Algebra]{On the Spectrum of Exterior Algebra, and Generalized Exponents of Small Representations}
\author[Sabino Di Trani]{Sabino Di Trani}
\address{Dipartimento di Matematica ``Guido Castelnuovo'', Sapienza - Universit\`a di Roma.
}
\email{sabino.ditrani@uniroma1.com }
\address{\emph{The autor has been partially supported by GNSAGA - INDAM group.}}
\address{ORCID id: \url{https://orcid.org/0000-0002-6651-558X}}

\begin{document}

\maketitle
\textbf{Abstract:}
We present some results about the irreducible representations appearing in the exterior algebra $\Lambda \g$, where  $\g$ is a simple Lie algebra over $\C$. For Lie algebras of type $B$, $C$ or $D$ we prove that certain irreducible representations, associated to weights characterized in a combinatorial way, appear as irreducible components of $\Lambda \g$. Moreover, we propose an analogue of a conjecture of Kostant, about irreducibles appearing in the exterior algebra of the little adjoint representation. Finally, we give some closed expressions, in type $B$, $C$ and $D$, for generalized exponents  of small representations that are fundamental representations and we propose a generalization of some results of De Concini, M\"{o}seneder Frajria, Procesi and Papi about the module of special covariants of adjoint and little adjoint type.  \\

\section{Introduction}
Let $\g$ be a simple Lie algebra over $\C$ of rank $\mathrm{rk} \,\g$. Fix a Cartan subalgebra $\h$ and let $\Phi$ be the associated root system with Weyl group $W$. We choose a  set of positive roots $\Phi^+$ and let $\Delta$ be the associated simple system. Let $\rho$ be the corresponding Weyl vector and $\theta$ the highest root with respect to the standard partial order $\leq$ on $\Phi^+$. If $\g$ is not simply laced, $\theta_s$ is the short dominant root. We denote by $\Pi$ and $\Pi^+$ the set of weights and the set of dominant weights respectively,  moreover we denote by $\omega_i$ the $i$-th fundamental weight. Throughout the paper, $V_\lambda$ will be the irreducible finite dimensional representation of $\g$ of highest weight  $\lambda \in \Pi^+$ and we denote by $V_\lambda^0$ the corresponding $W$-representation on the zero weight space of $V_\lambda$. Finally, $e_1\leq  \dots \leq e_n$ will be the \emph{exponents} of $\g$.\\
 The adjoint action of $\g$ on itself induces an action of $\g$ on $S(\g)$ and $\Lambda \g$, the symmetric and exterior algebras over $\g$ respectively, preserving the natural gradings.
Two celebrated results give an explicit description of the ring of invariants in $S(\g)$ and $\Lambda  \g$ with respect to this action.
\begin{theorem*}[Chevalley, Shephard and Todd]\label{thm:CST} Let $\g$ be a complex semisimple Lie algebra of rank $n$ and $\h$ a fixed Cartan subalgebra.
Up to identify $\g$ with $\g^*$ and $\h$ with $\h^*$ via Killing form, the restriction of polynomial functions induces an algebra isomorphism between the rigs of invariants
\begin{equation*}
 S(\g)^\g  \simeq S(\h)^W.
\end{equation*}
In particular,  $S(\g)^\g$ is a polynomial algebra with generators of degrees $e_1+1, \dots, e_n+1$.
\end{theorem*}
 \begin{theorem*}[Hopf, Koszul and Samelson]\label{thm:HKS}
Let $\g$ be a complex semisimple Lie algebra of rank $n$. Then
\begin{equation*}
\left( \Lambda \g \right) ^{\g} = \Lambda (P_1, \dots , P_{n}),
\end{equation*}
where the degree of a generator $P_i$ of the algebra of the invariants is equal to $2 e_i -1$.
 \end{theorem*}
If $M= \oplus M_i$ is a graded $\g$-module, we denote by \[P(V_\lambda, M,  t)= \sum_{i} \dim \Hom_\g (V_\lambda, M_i )t^i \]
the generating function for graded multiplicities of the irreducible representation $V_\lambda$ in $M$. As an immediate consequence of the above theorems, it is possible to obtain the following formulae that encode the graded structure of rings of invariants:
     \[P(V_0, \Lambda \g,  t)= \prod_{i=1}^n (1+t^{2e_i+1}), \qquad P(V_0, S(\g),  t)= \prod_{i=1}^n (1-t^{e_i+1})^{-1}.\]
Aiming to generalize the above results, some questions about irreducile representations in $S(\g)$ and $\Lambda \g$ naturally arise:
\begin{itemize}
 \item[\emph{Q1:}] Is it possible to determine the irreducible representations appearing in $S(\g)$ and in~$\Lambda \g$?
 \item[\emph{Q2:}] If $V_\lambda$ is a subrepresentation of $S(\g)$ or of $\Lambda \g$, is it possible to determine the degrees in which $V_\lambda$ appears?
  \item[\emph{Q3:}] Denoting by $\Lambda^i\g$ (resp. $S^i( \g)$) the submodule of homogeneous elements of degree $i$ in~$\Lambda \g$ (resp. $S(\g)$), is it possible to determine the multiplicity of $V_\lambda$ in $\Lambda^i \g$ (resp.~$S^i( \g)$)?
\end{itemize}
These questions inspired a great amount of claims and conjectures; many of them are still open or have only implicit answers. \\
For what concerns the irreducibles appearing in the symmetric algebra, the problem was extensively studied by Kostant in \cite{KostPoly}. More precisely, Kostant proved the isomorphism 
\[S(\g) \simeq S(\g)^\g \otimes \mathcal{H},\] where $\mathcal{H}$ is the ring of $\g$-harmonic polynomials, i.e. the ring of polynomials over $\g$ annihilated by $\g$-invariant differential operators of positive degreee with constant coefficients.
Studying the graded multiplicities of~$V_\lambda$ in $S(\g)$ can be then reduced to determining the multiplicity of $V_\lambda$ in each homogeneous component~$\mathcal{H}^i$ of~$\mathcal{H}$. 
Kostant proved that the multiplicity of $V_\lambda$ in $\mathcal{H}$ equals the dimension of~$V_\lambda^0$ and that the degrees $i$ such that~$V_\lambda$ appears in~$\mathcal{H}^i$ are related to the eigenvalues of the action of the Coxeter-Killing transformation on the $W$-representation~$V_\lambda^0$. \\
These integers are called the \emph{Generalized Exponents} associated to $V_\lambda$ and are extensively studied in the literature because of their nice combinatorial properties. We summarize some remarkable results about generalized exponents in Section \ref{sec:exp}.\\
On the other hand, despite its finite dimensionality, determining the irreducible components appearing in $\Lambda \g$ seems to be quite difficult. A complete description of irreducible representations in the exterior algebra is known only in type $A$, by some general arguments due to Berenstein and Zelevinsky, and for exceptional algebras of type $F_4$ and $G_2$, by direct computations. For other cases an open conjecture has been formulated by Kostant, describing the $V_\lambda$ appearing in $\Lambda \g$ as the irreducibles indexed by $\lambda$ smaller or equal to $2 \rho$ in the dominance order on weights, i.e. the ordering defined by the relation $\mu \leq \lambda$ if and only if $\lambda - \mu$ is a sum of positive roots. 
To introduce the reader to this fascinating subject and to provide a framework for the new results contained in this article, we present in Section \ref{sec:survey} a brief survey of some known results on this topic.\\
The remaining part of the paper is devoted to present our results.\\
In Section \ref{sec:BZ} we recall some results of Berenstein and Zelevinsky about multiplicities in tensor product decomposition. These techniques are used in \cite{BZ} to prove Kostant Conjecture in type~$A$. We use these tools to prove that large families of irreducible representations appear as irreducible representations in $\Lambda \g$, for $\g$ of type $B$, $C$ and $D$. More precisely we introduce the Coordinatewise Ordering (Definition \ref{defn:lexorder}) on the set of dominant weights, prescribing that $\mu$ is \emph{coordinatewise smaller} than $\lambda$ (for short $\mu \lesssim \lambda$) if certain combinatorial conditions are satisfied. We use this ordering to describe a suitable subset of the set of dominant weights smaller than $2 \rho$ in the dominance order. We prove that irreducible representations associated to weights in this subset appear in $\Lambda \g$.  The main result of the section is the following theorem:
\begin{theorem*}Let $\g$ be a simple Lie algebra over $\C$ of type $B$, $C$ or $D$ and let $\lambda$ be a dominant weight for $\g$. If $\lambda \leq 2 \rho$ and $\lambda \lesssim 2 \rho$, then $V_\lambda$ appears as irreducible component in $\Lambda \g$.\end{theorem*}
Section \ref{sec:small} is devoted to compute explicit formulae for polynomials of generalized exponents, using the techniques summarized in Section~\ref{sec:exp}.
In particular, denoting by $E_\lambda(t)$ the generating polynomial of generalized exponents associated to $V_\lambda$, i.e. the Poincarè polynomial of graded multiplicities of $V_\lambda$ into $\mathcal{H}$, in Section \ref{sec:small} we observe that the following formula can be obtained in type $C_n$ as a consequence of results contained in \cite{SDT}
\[ E_{\omega_{2k}}(t)=\frac{t^{2k}(n-2k+1)_{t^2}}{(n-k+1)_{t^2}}\binom{n}{k}_{t^2}, \]
where $(n)_t$ denotes the $t$-analogue of $n$ and $\binom{n}{k}_t$ is the $t$-binomial.
Moreover, denoting by $\lfloor a \rfloor$ the integer part of $a$, we prove that the following formulae hold in type $B_n$
\[E_{\omega_{2k}}(t)=t^k\binom{n}{k}_{t^2}, \]\[ E_{\omega_{2k+1}}(t)=t^{n-k}\binom{n}{k}_{t^2},\]\[ E_{2\omega_{n}}(t)=t^{n-\lfloor\frac{n}{2}\rfloor}\binom{n}{\lfloor\frac{n}{2}\rfloor}_{t^2}, \]
and in type $D_n$
\[E_{\omega_{2k}}(t)=t^k\frac{(t^{n-2k}+1)}{(t^n+1)}\binom{n}{k}_{t^2},\]
\[ E_{\omega_{n-1}+\omega_n}(t)=\frac{t^{\lfloor\frac{n}{2}\rfloor}(t+1)}{(t^n+1)}\binom{n}{\lfloor\frac{n}{2} \rfloor}_{t^2}, \]\[ E_{2\omega_{n-1}}(t)=E_{2\omega_{n}}(t)=\frac{t^{\frac{n}{2}}}{(t^n+1)}\binom{n}{\frac{n}{2}}_{t^2},\]
where the formula for $E_{\omega_{n-1}+\omega_n}(t)$  holds for $n$ odd and the formulae for $E_{2\omega_{n-1}}(t)$ and $E_{2\omega_{n}}(t)$ must be considered only if $n$ is even. Finally, some open question and conjectures are proposed at the end of Sections \ref{sec:BZ} and of Section \ref{sec:small}.
\subsection*{Acknowledgements} The main original contributions of this paper are some results that I obtained during  my doctoral studies, so I would like to thank my advisor, Professor Paolo Papi, for his mentoring and for his supervision. Moreover, I am grateful to Professor Andrea Maffei for many useful discussions about the Kostant Conjecture. I would like to extend my special thanks to the anonymous referee for their really careful reading and for their precious comments to a previous version of the paper. I am also grateful to Rosario Mennuni and Viola Siconolfi for
their advice on the organization of a first draft of the paper.  
Finally, this article was partially written during my frequent stays in Pisa: I express my gratitude to P.F., to M.A.P. and to the little R.F. for their great hospitality and to all my fiends at the Mathematics Department for their support.
\section{Irreducible Representations in the Exterior Algebra}\label{sec:survey}
As mentioned in the introduction, an uniform description of irreducible representations appearing in the exterior algebra $\Lambda \g$, with $\g$ a simple Lie algebra over $\C$, has been proposed by Kostant:
 \begin{conj}[Kostant, c.f.r.~\cite{BZ1}, Introduction]\label{conj:Kost}
  The representation $V_\lambda$ appears in the decomposition of $\Lambda \g$ if and only if $\lambda \leq 2\rho$ in the dominance order.
 \end{conj}
Currently a proof of this conjecture is known only in type $A$, by the combinatorial construction given in \cite{BZ1}, and in the exceptional cases $G_2$ and $F_4$ by explicit computations, as reported in \cite{CKM}. Moreover, we mention that in \cite{CKM} the authors exhibit a possible uniform proof of the Kostant Conjecture for algebras of types $ADE$, assuming that 1 is a saturation factor for any simply laced algebra. It is not clear if similar techniques could be used to prove the Conjecture in the remaining cases.  
Moreover, \emph{a priori} it should be possible to verify Kostant Conjecture in type $E$ by direct computations, but it seems to be an unfruitful approach. Nevertheless, a uniform proof of Conjecture \ref{conj:Kost} is desirable, but a concrete strategy is far to be clear. In addition to that, if $V_\lambda$ appears in $\Lambda \g$ studying its graded multiplicities seems to be also very complex. We collect here some partial related results. Firstly, a uniform bound for multiplicity of $V_\lambda$ is known.
 \begin{thm}[Reeder, \cite{Reeder}, Section 4]\label{thm:multsmall}
\begin{equation}\label{multReeder} \dim \,\mathrm{Hom}_{\g}(V_\lambda, \Lambda \g) \leq 2^{\mathrm{rk}\g} \dim V_\lambda^0.\end{equation}
\end{thm}
Moreover, Reeder investigated when the equality holds.
\begin{defn}[c.f.r \cite{Reeder}, Definition 2.2]
An irreducible representation $V_\lambda$ is \emph{small} if $\lambda$ is in the root lattice and if $2 \alpha \nleq \lambda $ for every dominant root $\alpha$.
\end{defn}
 \begin{thm}[Reeder, \cite{Reeder}, Section 4]
Equality in Formula (\ref{multReeder}) holds if and only if $V_\lambda$ is small.
\end{thm}
Observe in particular that the adjoint and the little adjoint representations are special cases of small representations.
Some explicit formulae for polynomials of graded multiplicities are proved by Bazlov.
\begin{thm}[Bazlov, \cite{Baz}, Section 5.2] The following formula for graded multiplicities of adjoint representation in $\Lambda \g$ holds:
\[P(\g, \Lambda \g, q)= (1+q^{-1})\prod_{i=1}^{n-1}(q^{2e_i+1}+1) \sum_{i=1}^n q^{2 e_i}. \]
\end{thm}
Moreover, for certain weights close to $2 \rho$, an explicit formula 
can be found in \cite{Reeder}.
\begin{thm}[Reeder, \cite{Reeder}, Proposition 6.3]\label{thm:gradmultdI}
Let $ I \subseteq \Delta $. Consider $\delta_I = \sum_{\alpha \in I} \alpha$ and denote by $c(I)$ the number of connected component of the Dynkin subdiagram generated by $I$. Then
\[P(V_{2\rho -\delta_I}, \Lambda \g, t)=t^{|\Phi^+|-|I|}(t+1)^{n-c(I)}(t^2+1)^{|I|-c(I)}(t^3+1)^{c(I)}.\]
\end{thm}
Similarly, closed formulae can be obtained for small representations as a consequence of a conjecture formulated by Reeder in \cite{Reeder} and proved in \cite{SDT} and \cite{SDT2}. This conjecture was inspired by two remarkable results:
\begin{theorem}[Broer \cite{Broer}, Theorem 1]\label{Broer}
The homomorphism induced by the Chevalley restriction theorem
\[\mathrm{Hom}_\g(V_\lambda, S(\g)) \rightarrow \mathrm{Hom}_W(V^0_\lambda, S(\h))\]
is a graded isomorphism of $S(\g)^\g \simeq S(\h)^W$-modules if and only if $V_\lambda$ is \emph{small}.
\end{theorem}
\begin{theorem}[Chevalley, Eilenberg \cite{CE}, Reeder \cite{R1}]\label{thm:cohomology}
 Let $G$ be a compact Lie group, $T \subset G$ a maximal torus and $W$ its Weyl group. Let $\g$ be the complexified Lie algebra of $G$ and $\h$ the Cartan subalgebra of $\g$ associated to $T$. The Weyl map $\psi: G/T\times T \rightarrow G$ induces in cohomology the following graded isomorphism:
 \begin{equation*}
  \left(\Lambda \g\right)^\g \simeq H^*(G) \simeq \left(H^*(G/T) \otimes H^*(T)\right)^W \simeq \left(  \mathcal{H}_{(2)} \otimes \Lambda \h \right)^W.
 \end{equation*}
 where $\mathcal{H}_{(2)}$ denotes the graded ring of $W$-harmonic polynomials over $\h$, with a grading obtained by doubling the natural one.
\end{theorem}
Theorem \ref{Broer} and Theorem \ref{thm:cohomology} suggest that graded multiplicities of a small representation $V_\lambda$ in $\Lambda \g$ are linked to multiplicities of the $W$-representation $V_\lambda^0$ in the bigraded ring $\Lambda \h \otimes \mathcal{H}_{(2)}$.
Reeder conjectured that, if $V_\lambda$ be a small representation, the following equality holds:
\begin{equation}\label{form:reeder}\dim \mathrm{Hom}_{\g} (V_\lambda, \Lambda^i \g )=  \sum_{k+h=i}\dim \mathrm{Hom}_{W} (V_\lambda^0,  \mathcal{H}^h_{(2)} \otimes \Lambda^k \h )\end{equation}
Small representations for algebras of type $A_{n-1}$ are of the form $V_\lambda$ where $\lambda$ is a partition of~$n$. Reeder's conjecture is implicitly proved in literature for algebras of type $A$ by comparing the results contained in \cite{KP} and \cite{M} with the following formula proved by Stembridge:
\begin{thm}[Stembridge, \cite{St1}, Corollary 6.2] Let $\lambda$ be a partition of $n$ and $\Gamma$ the associated Young tableaux, displayed in the English way.
  \[P(V_\lambda , \Lambda \g , q) = \frac{\prod_{i=1}^n (1-q^{2i})}{(1+q)} \prod_{(ij) \in \Gamma }\frac{\left( q^{2j-2}+q^{2i-1}\right)}{\left(1-q^{2h(ij)}\right)}\]
  where $h(ij)$ denotes the hook length of the box $(ij) \in \Gamma$.
 \end{thm}
For other simple Lie algebras the conjecture is proved in \cite{SDT} and \cite{SDT2} using a case by case strategy. The problem of finding a uniform approach to prove Equation~(\ref{form:reeder}) for small representations is still open and very interesting.
In this spirit, an enhanced version of Reeder's conjecture has been recently proposed in \cite{DCP}, Section 7.\\
Finally, we remark that the module of special coinvariants $\Hom_\g (\g, \Lambda \g)$ has a richer geometric structure, as proved in~\cite{DCPP}:
 \begin{thm}[De Concini, Papi, Procesi, \cite{DCPP}, Theorem 1.1]\label{thm:DCPP}
The module $Hom_\g (\g , \Lambda \g)$ is a finitely generated free module over $\Lambda(P_1, \dots, P_{n-1})$ with generators in degree $2e_i$ and $2e_i-1$.
 \end{thm}
An analogous result it is proved in~\cite{DCMPP}, when $\g$ is not simply laced, for the module $\Hom_\g (V_{\theta_s}, \Lambda \g)$.
An extension of these theorems to certain small representations is proposed in Section~\ref{subsec:conj}.
\section{Berenstein and Zelevinsky Polytopes}\label{sec:BZ}
The more efficient way to approach the Kostant Conjecture seems to be by facing the problem using tensor product decomposition techniques. In fact, using the Weyl Character Formula, in \cite{Kost} Kostant proved the following isomorphism:
\[\Lambda \g \simeq \left(V_\rho \otimes V_\rho \right)^{\oplus 2^{\mathrm{rk} \g}} \]
Kostant's Conjecture can be consequently reformulated in the following terms (c.f.r.~\cite{CKM}, Remark 4):
 \begin{conj}[Kostant]
  The representation $V_\lambda$ appears in the decomposition of $V_\rho \otimes V_\rho$ if and only if $\lambda \leq 2\rho$ in the dominance order.
 \end{conj}
The conjecture in type $A$ is proved by Berenstein and Zelevinsky in \cite{BZ1} as a consequence of a more general combinatorial construction, used to find the tensor product decomposition of two irreducible finite dimensional representations of $\mf{gl}_n (\C)$. More in detail, they prove that for any triple of dominant weights $(\lambda, \mu, \nu)$, the irreducible representation $V_\nu$ is a component of $V_\lambda \otimes V_\mu$ if and only if there exists an integral point in a suitable polytope $P(\lambda, \mu , \nu)$ depending on the expansion of $\lambda$ and $\mu$ in terms of the fundamental weights.
As an application of their results, Berenstein and Zelevinsky prove that for every $\mu \leq 2 \rho$ the polytopes of the form $P(\rho, \rho , \mu)$ have \emph{at least} one integral point.
Moreover in \cite{BZ} it is conjectured that a similar description of tensor multiplicities in terms of integral points of certain polytopes holds for every classical Lie algebra. The statement of the conjecture is recalled in subsection \ref{subsec:BZpart}, it is proved by Berenstein and Zelevinsky as a consequence of results contained in \cite{BZ2}.
\subsection{Orderings on Dominant Weights.}\label{subsec:weights}
We recall now how roots systems of type $B_n$, $C_n$ and $D_n$ can be realized in an $n$-dimensional euclidean vector space $\mathbb{E}$ with basis $\{\varepsilon_1, \dots, \varepsilon_n\}$. We follow the constructions exposed in \cite{Bou} and \cite{FH}.\\
\emph{Root System of Type $B_n$:}
\[\Phi=\{\pm \varepsilon_i\pm \varepsilon_j\}_{i<j} \cup \{\pm \varepsilon_1, \, \dots \, , \pm \varepsilon_n\},\]
\[\Delta=\{\varepsilon_1-\varepsilon_2, \,\dots, \,  \varepsilon_{n-1}-\varepsilon_n, \,\varepsilon_n\}, \]
\[\Phi^+=\{ \varepsilon_i\pm \varepsilon_j\}_{i<j} \cup \{ \varepsilon_1, \, \dots \, ,  \varepsilon_n\}\quad W=S_n\ltimes \left(\Z/2\Z\right)^n,\]
\[\omega_i= \varepsilon_1+ \dots +\varepsilon_i \quad  \omega_n= \frac{\varepsilon_1+ \dots +\varepsilon_n}{2},\]
\[ \rho = \frac{(2n-1)\varepsilon_1 + (2n-3)\varepsilon_2 + \dots + 3\varepsilon_{n-1}+\varepsilon_n}{2}. \]
\emph{Root System of Type $C_n$:}
\[\Phi=\{\pm \varepsilon_i\pm \varepsilon_j\}_{i<j} \cup \{\pm 2\varepsilon_1, \, \dots \, , \pm 2\varepsilon_n\},\]
\[\Delta=\{\varepsilon_1-\varepsilon_2, \,\dots, \,  \varepsilon_{n-1}-\varepsilon_n, \,2\varepsilon_n\}, \]
\[\Phi^+=\{ \varepsilon_i\pm \varepsilon_j\}_{i<j} \cup \{ 2\varepsilon_1, \, \dots \, ,  2\varepsilon_n\}\quad W=S_n\rtimes \left(\Z/2\Z\right)^n,\]
 \[\omega_i=\varepsilon_1+ \dots + \varepsilon_i,\]
 \[ \rho = n\varepsilon_1 + (n-1)\varepsilon_2 + \dots + 2\varepsilon_{n-1}+\varepsilon_n. \]
\emph{Root System of Type $D_n$:}
\[\Phi=\{\pm \varepsilon_i\pm \varepsilon_j\}_{i<j} \quad \Delta=\{\varepsilon_1-\varepsilon_2, \,\dots, \,  \varepsilon_{n-1}-\varepsilon_n, \,\varepsilon_{n-1}+\varepsilon_n\}, \]
\[\Phi^+=\{ \varepsilon_i\pm \varepsilon_j\}_{i<j} \quad {W=S_n\ltimes \left(\Z/2\Z\right)^{n-1}},\]
 \[\omega_i= \varepsilon_1+ \dots +\varepsilon_i\quad \omega_{n-1}= \frac{\varepsilon_1+ \dots -\varepsilon_n}{2} \quad \omega_n= \frac{\varepsilon_1+ \dots +\varepsilon_n}{2},\]
 \[\rho = (n-1)\varepsilon_1 + (n-2)\varepsilon_2 + \dots + \varepsilon_{n-1}.\]
The set dominant weights is partially ordered by the dominance order, i.e. $\lambda \geq \mu$ if $\lambda - \mu$ is a sum of positive roots.
Moreover, every dominant weight $\lambda$ can be written as a sum $\lambda_1 \varepsilon_1 + \dots \lambda_n \varepsilon_n$ where $\lambda_i \in \frac{1}{2} \Z$ for all $i$.
The condition $\lambda \geq \mu$ in the dominance order can be restated as follows:
\begin{rmk}\label{classmin2rho}
Let $\lambda=\lambda_1 \varepsilon_1 + \dots +\lambda_n \varepsilon_n$ and $\mu=\mu_1 \varepsilon_1 + \dots +\mu_n \varepsilon_n$ be two dominant weights for a simple Lie algebra of type $B_n$, $C_n$ or $D_n$. Then $\lambda \geq \mu$ if and only if the following conditions hold:
\begin{enumerate}
\item $\sum_{i=1}^k (\lambda_i-\mu_i) \geq 0 $ for all $1 \leq k \leq n$, in type $B$;
 \item $\sum_{i=1}^k (\lambda_i-\mu_i) \geq 0 $ for all $1 \leq k \leq n$
and  $\sum_{i=1}^n (\lambda_i - \mu_i )$ is an even integer, in type $C$ and $D$.
\end{enumerate}
\end{rmk}
We introduce now a different ordering on the set of weight.
\begin{defn}[Coordinatewise order on weights]\label{defn:lexorder}
 Let $\lambda=\lambda_1 \varepsilon_1 + \dots +\lambda_n \varepsilon_n$ and $\mu=\mu_1 \varepsilon_1 + \dots +\mu_n \varepsilon_n$ be two dominant weights for a simple Lie algebra $\g$ of type $B_n$,$C_n$ or $D_n$. We say that $\mu$ is smaller than $\lambda$ with respect to the relation $\lesssim$ if and only if $\lambda_i - \mu_i \geq 0$ and $|\lambda_i|\geq|\mu_i|$ for all $i$. In this case we write $\mu \lesssim \lambda$ and we say that $\mu$ is smaller than $\lambda$ with respect to the coordinatewise order.
\end{defn}
\begin{rmk}
Observe that the coordinatewise ordering is different from the dominance ordering. As an example, in type $C$ the weight $\omega_2$ is the only non zero dominant weight smaller than $2 \omega_1$ with respect to the dominance order, but $\omega_2 \not\lesssim 2\omega_1$. On the other side in type $C$ we have that $\omega_1 \lesssim \omega_2$, although $\omega_2$ is a minimal element between non zero dominant weights with respect to dominance order. 
\end{rmk}
The next two sections are devoted to prove the following theorem:
\begin{thm}\label{thm:lexdom}Let $\g$ be a simple Lie algebra over $\C$ of type $B$, $C$ or $D$ and let $\lambda$ be a dominant weight for $\g$. If $\lambda \leq 2 \rho$ and $\lambda \lesssim 2 \rho$, then $V_\lambda$ appears as irreducible component in $\Lambda \g$.\end{thm}
\begin{example} In this example we compare the set of weights considered in the assert of Theorem~\ref{thm:lexdom} with the ones appearing in Theorem~\ref{thm:multsmall} and in Theorem~\ref{thm:gradmultdI}. In particular we focus on the case of simple Lie algebra $C_3$. In type $C_3$ there are 35 dominant weights smaller or equal to $2 \rho$ with respect to dominance order. Between these weights, there are 30 dominant weights $\mu$ such that $\mu \lesssim 2 \rho$. All small weights appear in this set, but they are considerably fewer (more precisely, in type $C_3$ there are 4 small weights, c.f.r. Table 4). Moreover, in type $C_3$ there are 7 dominant weights of the form $2\rho - \delta_{I}$ with $I \subset \Delta$. Between them only 4 weights are not smaller than $2 \rho$ with respect to the coordinatewise order. 
\end{example}
\subsection{$\g$-partitions and Berenstein-Zelevinsky polytopes.}\label{subsec:BZpart}
Let $m$ be a weight in the root lattice for the Lie algebra $\mf{so}_{2n+1} \C$, it can be described by a vector of non negative integers \[(m_{12}, m^+_{12}, \dots, m_{n-1n},m^+_{n-1n},m_1, \dots,m_n)\] such that
  \[m=\sum_{i<j}m_{ij}(\varepsilon_i-\varepsilon_j)+\sum_{i<j}m^+_{ij}(\varepsilon_i+\varepsilon_j)+\sum_{i}m_i\varepsilon_i.\]
  We say that the sequence of integers $(m_{12}, m^+_{12}, \dots, m_{n-1n},m^+_{n-1n},m_1,\dots,m_n)$ is an \emph{$\mf{so}_{2n+1}$-partition for $m$}.
  We say that an $\mf{so}_{2n+1}$ partition is an $\mf{so}_{2n}$-partition (resp. $\mf{sp}_{2n}$-partition) if $m_i=0$ (resp. $m_i$ is even) for every $i$.
The inequalities that determine the Berenstein-Zelevinsky polytope for a general tensor product $V_\lambda \otimes V_\mu$ are described in  \cite{BZ} in terms of the variables $m_{12}, m^+_{12}, \dots, m_{n-1n}, m^+_{n-1n}$ and $m_1, \dots, m_n$.  We recall here their  description as presented in \cite{BZ}.\\
Consider the set $I=\{\bar{0}, 1, \dots , n, \bar{1}, \dots, \bar{n}\}$, ordered by $\bar{0}<1<\bar{1}< \dots < n < \bar{n}$, and set
\begin{equation*}
\Delta_{ij}=m_{ij}-m^+_{ij}, \quad \Delta_{\bar{i} \bar{j}}=\Delta_{i+1j+1},
\quad 
 \Delta_{i \bar{j}}=\Delta_{\bar{i}j}=
\left\{
	\begin{array}{lll}
	m^+_{i,j+1}-m_{i+1 j+1} & \mbox{if } j<n, \\
m_i-m_{i+1}  & \mbox{if } i=n.\\
	\end{array}
\right.
\end{equation*}
where $m_{i,j},m^+_{i,j}$ must be considered only if $i<j$.
Now, if $j<n$ and $t \in I$, we consider the linear forms (c.f.r.~\cite{BZ}, Formulae (2.4)):
\begin{equation}\label{j<n}
\left \{
\begin{array}{ll}
\mathscr{L}^t_j(m)=-\sum_{\bar{0} \leq s \leq t} \Delta_{sj},   \\
\mathscr{N}^{t,0}_j(m)= \Delta_{\bar{j}j}+ \sum_{j+1 \leq s \leq t} \Delta_{\bar{j},s},   \\
\mathscr{N}^{t,1}_j(m)= \mathcal{N}^{n,0}_j + \sum_{t \leq s \leq n}  \Delta_{j,s}.   \\
	\end{array}
\right.
\end{equation}
Otherwise, if $j=n$, consider 
\begin{equation}\label{eqTB}
     \mathscr{L}_n^{t}(m)=-\left[   2 \left(    \sum_{1 \leq p \leq t } \Delta_{p\,n}\right)  +    \sum_{0 \leq p \leq t } \Delta_{\overline{p}\,n} \right]  \qquad \mathcal{N}^{n,1}_n(m)= m_n \qquad \mbox{(Type B)},
    \end{equation}
\begin{equation}\label{eqTC}
     \mathscr{L}_n^{t}(m)=-\left[    \left(    \sum_{1 \leq p \leq t } \Delta_{p\,n}\right)  +  \left(  \frac{1}{2}  \sum_{0 \leq p \leq t } \Delta_{\overline{p}\,n}\right) \right] \qquad \mathcal{N}^{n,1}_n(m)= m_n/2 \qquad \mbox{(Type C)},
    \end{equation}
    \begin{equation}\label{eqTD}
 \mathscr{L}_n^{t}(m)= {\widehat{\mathscr{L}}_{n-1}}^{t}(m) \qquad \mathcal{N}^{n,1}_n(m)= m^+_{n-1,n} \qquad \mbox{(Type D)},
\end{equation} 
where $\widehat{\mathscr{L}}_{n-1}^{t}(m)$ is the image of $\mathscr{L}_{n-1}^{t}(m)$ under the involution 
\begin{equation*}
\widehat{m}_{i,j}=
\left \{
\begin{array}{ll}
m_{i,j} \quad \mbox{if $j<n$}  \\
m^+_{i,j} \quad  \mbox{if $j=n$}  \\
\end{array}
\right.
\qquad 
{\widehat{m}^+}_{i,j}=
\left \{
\begin{array}{ll}
m^+_{i,j} \quad \mbox{if $j<n$}  \\
{m}_{i,j} \quad \mbox{if $j=n$}  \\
\end{array}
\right.
\end{equation*}
Let us denote with $c_{\lambda \mu}^\nu$ the generalized Littlewood Richardson coefficient associated to the triple of dominant weights $( \lambda, \mu, \nu)$, i.e. the multiplicity of $V_\nu$ in $V_\lambda \otimes V_\mu$. The following theorem, crucial for our results, was conjectured in \cite{BZ} and proved in \cite{BZ2}
 \begin{theorem}[Berestein, Zelevinsky, \cite{BZ}, Section 2]\label{BZ}
Let $\lambda=a_1 \omega_1 + \dots + a_n \lambda_n$  and $\mu = b_1 \omega_1 + \dots + b_n \omega_n$ be dominant weights.
The irreducible components of $V_\lambda \otimes V_\mu$ are in bijection with integral points of the polytope defined by the inequalities
\begin{equation*}
\mathscr{L}^t_j \leq a_j  \qquad \mathscr{N}^{t, 0}_j  \leq b_j \qquad \mathscr{N}^{t, 1}_j \leq b_j, 
\end{equation*}
where the indices considered are displayed in the Table~\ref{tab:indices}.
\end{theorem}
\begin{table}[ht!]
  \begin{center}
   \caption{Indices contribution}\label{tab:indices}
   \begin{tabular}{|c|c|c|}
   \hline
 & \textbf{Type B and C}  & \textbf{Type D} \\
    \hline
\multirow{2}{*}{$\mathscr{L}^t_j$}&\multirow{2}{*}{$1 \leq j \leq n, \; \bar{0} \leq t <j$}&\multirow{1.2}{*}{$1 \leq j \leq n-1, \; \bar{0} \leq t <j$} \\ & & $j=n, \bar{0} \leq t <n-1$\\
\hline
\multirow{2}{*}{$\mathscr{N}^{t,0}_j$}&\multirow{2}{*}{$1 \leq j \leq n-1, \, \bar{j} \leq t \leq n$}&\multirow{2}{*}{$1 \leq j \leq n-2, \, \bar{j} \leq t \leq n-1$}\\
& &\\
\hline
\multirow{2}{*}{$\mathscr{N}^{t,1}_j$}&\multirow{1.2}{*}{$1 \leq j \leq n-1, \, \bar{j} < t \leq n$}  &\multirow{1.2}{*}{$1 \leq j \leq n-2, \, \bar{j} < t \leq n$,} \,  \\ &\multirow{1}{*}{ $j=t=n$} & $j=t=n$, \; $j=n-1, \, t=n$\\
\hline
   \end{tabular}
\end{center}
\end{table}
Each integral point in the polytope corresponds to a $\g$~partition. We are going to call these $\g$~partitions  \emph{admissible for the pair $(\lambda, \mu)$}. We say that a $\g$~partition $(m_{12}, m^+_{12}, \dots,m_n)$ is associated to a weight $\nu$ if
  \[\nu=\sum_{i<j}m_{ij}(\varepsilon_i-\varepsilon_j)+\sum_{i<j}m^+_{ij}(\varepsilon_i+\varepsilon_j)+\sum_{i}m_i\varepsilon_i\]
As a corollary of the Theorem~\ref{BZ}, Berenstein and Zelevinsky prove that:
\begin{theorem}[Berestein, Zelevinsky, \cite{BZ}, Section 2]
The coefficient $c_{\lambda \mu}^\nu$ is equal to the number of $\g$-partitions admissible for the pair $(\lambda, \mu)$ and associated to $\lambda + \mu - \nu$.
\end{theorem}
We want to use the previous results to obtain informations about the decomposition into irreducibles of $V_\rho \otimes V_\rho$.
In particular, studying the irreducible components which appear in $V_\rho \otimes V_\rho$ is consequently equivalent to describe the integral points in the polytope defined by
\begin{equation*}\label{BZeq}
\mathscr{L}^t_j(m)  \leq 1,  \qquad \mathscr{N}^{t, 0}_j(m), \,    \mathscr{N}^{t, 1}_j(m) \leq 1,
\end{equation*}
for $t,j$ that range as in Table~\ref{tab:indices}. From now on this section, by abuse of notation, we say that a $\g$~partition is \emph{admissible} if it is admissible for the pair $(\rho, \rho)$.
Our aim is to construct explicitly an admissible $\g$-partition associated to each weight $\lambda \leq 2 \rho$ such that $\lambda \lesssim 2 \rho$.\\
Firstly, we rearrange the equations defining the Berenstein and Zelevinsky polytopes in a more explicit form.
Set $M(i,j)=m_{ij}-m^+_{ij}$, $N(i)=m_i-m_{i+1}$, $R(i,j)=m^+_{i,j}-m^+_{i+1,j}$ and $S(i,j)=m_{ij+1}-m_{i+1 \, j+1} + m^+_{ij+1}-m^+_{i+1 \, j+1}$ for $j\in\{1, \dots, n\}$ and $1\leq i < j$, then the linear forms in Formula~(\ref{j<n}) can be expressed as:
\[\mathscr{L}_j^t(m)=    \sum_{i=1}^{t-1 } \left( M(i,j+1) - M(i, \, j)\right) - M(t,j)   + m_{t \, j+1},\]
  \[\mathscr{L}_j^{\overline{t}}(m)=    \sum_{i=1}^{t } \left( M(i,j+1) - M(i, \, j)\right)    + m_{t +1\, j+1} ,\]
   \[\mathscr{N}_i^{t \; 0}(m) = m^+_{i\, i+1} + \sum_{j=i +1  }^{t-1} R(i,j+1) + (m^+_{i \, t+1} - m_{i+1 \, t+1}),\]
   \[\mathscr{N}_i^{\overline{t} \; 0}(m) = m^+_{i\, i+1} + \sum_{j=i +1  }^{t} R(i,j+1),\]
   \[\mathscr{N}_i^{n \; 0} (m)= m^+_{i\, i+1} + \sum_{j=i +1  }^{n-1} R(i,j+1) + N(i) ,\]
   \[\mathscr{N}_i^{t \; 1}(m)= m^+_{i\, i+1} + N(i)  + M(i,t) + \sum_{j=i +1  }^{t-1} R(i,j+1) +\sum_{j=t}^{n-1} S(i,j+1) , \]
 \[\mathscr{N}_i^{\overline{t} \; 1}(m)=  m^+_{i\, i+1} + N(i) + \sum_{j=i +1  }^{t-1} R(i,j+1) +\sum_{j=t}^{n-1} S(i,j+1).  \]
If $j=n$, the equations (\ref{eqTB}), (\ref{eqTC}) and (\ref{eqTD}) can be rearranged in the following way:
\begin{equation*}\label{ReqTB}
   \mathscr{L}_n^t(m)=  - 2\sum_{i=1}^{t }  M(i, \, n)    + m_{t}  \qquad \mathscr{L}_n^{\overline{t}}(m)=  - 2\sum_{i=1}^{t }  M(i, \, n)    + m_{t+1} \qquad \mbox{(Type B)},
    \end{equation*}
\begin{equation*}\label{ReqTC}
     \mathscr{L}_n^{\overline{t}}(m)=  - \sum_{i=1}^{t }  M(i, \, n)   + m_{t+1}/2 \qquad \mathscr{L}_n^t(m)=  - \sum_{i=1}^{t }  M(i, \, n)    + m_{t}/2 \qquad  \mbox{(Type C)},
    \end{equation*}
\footnotesize{
\begin{equation*}\label{ReqTD}
\mathscr{L}_n^t(m)=-\sum_{i=1}^{t-1}M(i, \, n)-\sum_{i=1}^{t } M(i, \, n-1)+ m^+_{t \, n}\qquad \mathscr{L}_n^{\overline{t}}(m)=   - \sum_{i=1}^{t } \left( M(i, \, n) + M(i, \, n-1)\right)    + m^+_{t +1\, n} \quad \mbox{(Type D)}.
\end{equation*} 
}
\normalsize
Here we adopted the convention that, if the set of indices is empty, the sum is equal to 0.
\subsection{The construction.}\label{subsec:construction}
 For each $\lambda \leq 2 \rho$ set $c_i=2|\rho_i| - |\lambda_i|$, where by $\lambda_i$ and $\rho_i$ we denote the $i$-th coordinate of $\lambda$ and $\rho$, with respect to the basis $\{\varepsilon_1, \dots \varepsilon_n\}$. If $0 \leq c_i$ for all $i \leq n$, we give an explicit construction of an admissible $\g$-partition associated to $2\rho-\lambda$, appearing as integral point in the Berenstein Zelevinzky polytope associated to the tensor product $V_\rho \otimes V_\rho$. The conditions on the $c_i$ in particular are equivalent to require that $\lambda \lesssim   {2 \rho}$. \\
We have three main cases, depending on the parity of the $\{c_i\}_{i\leq n}$. We will construct an admissible $\g$-partition $m=(m_{12}, \dots , m_n)$ associated to $2 \rho - \lambda $ in an iterative way. We start setting $m$ to be the zero vector.\\
\textbf{Case A: the $c_i$ are all even.} 
\begin{itemize}
 \item[\emph{Step 1}]  If $c_n=0$ set $m_n=0$, otherwise $m_n=2$ (observe that the case $c_n$ even and greater than  $0$  cannot happen in type $B$ and $D$  because in these cases $2\rho_n<2$);
 \item[\emph{Step $h+1$}] Suppose $h+1=n-(i-1)+1$ and let  $(m_{i \, i+1}, m^+_{i \, i+1 }, \dots , m_{i \, n}, m^+_{i \, n }, m_i)$ be the integers constructed  at the $h$-th step. Let $J_i=\{j_k < \dots < j_1\}$ be the set of indices such that $m_{ij_s}\neq 0$. By convention, we set $j_0=n+1$. We have the following cases:
 \begin{enumerate}
 \item if $c_{i-1}=0$, set $m_{i-1 , j}=m^+_{i-1 , j}=m_{i-1}=0 $ for all $j$;
  \item if $c_i \geq c_{i-1}>0$, set $m_{i-1}=m_i$,  $m_{i-1 \, j}=m_{i \, j}$ and $m^+_{i-1 \, j}=m^+_{i \, j}$ for all $j$ such that $n \geq j \geq j_s$, where $s$ is chosen to be equal to $c_{i-1}/2 $ if  $m_i=0$, and to  $c_{i-1}/2 -1$ otherwise.
  Finally set $m_{i-1j}=m^+_{i-1j}=0$ for the remaining indices;
  \item if $c_{i-1}=c_i+2$, set $m_{i-1}=m_i$ and $m_{i-1 \, j}=m_{i \, j}$ e $m^+_{i-1 \, j}=m^+_{i \, j}$ for all $j > i$. Finally set $m_{i-1 \, i}=m^+_{i-1 \, i}=1$.
 \end{enumerate}
\end{itemize}
\begin{proposition} 
The construction exposed in Case A produces an admissible $\g$-partition associated to~$2 \rho  - \lambda$. \end{proposition}
\proof By Theorem~\ref{BZ}, we need to prove that $\mathscr{L}^i_j(m), \mathscr{N}^{i, 0}_j(m)$ and $\mathscr{N}^{i, 1}_j(m)$ are smaller than~1.  Observe that in our construction  $m_i \neq 0$ only if $m_{i+1}\neq 0$ and then $N(i) \leq 0$ for all $i$. Moreover, a non zero $m_{ij}$ is constructed (i.e in case~$(2)$ or in case~$(3)$) if and only if $m^+_{ij}\neq 0$, and in that case we always have $m_{ij}=m^+_{ij}$. Consequently $M(i,j)=0$ for every pair $i,j$. Finally, if $i +1 < j$, we always have that $m_{i,j}=m^+_{ij} \neq 0$ only if $m_{i+1,j} =m^+_{i+1j}\neq 0$ and then $R(i,j), \, S(i,j) \leq 0 $. Verify that the constructed  $\g$-partition is admissible is now just a straightforward computation, recalling that by construction described in $(2)$ and $(3)$ we have $m_{ij}, m^+_{ij}\leq 1$ for every pair $i,j$ and $m^+_{ij} -m_{i+1j} \leq 0 $ for every $j$ such that $i+1 <j$.
\endproof
\begin{example}\label{ex:gparteven}
In this example we construct admissible $\mf{sp}_6\C$-partitions associated to the weights $2 \rho - \lambda$ and $2 \rho - \lambda'$, where $\lambda=2 \omega_3$ and $\lambda' = 4 \omega_1 $. We remark that, because in type $C_3$ we have nine positive roots, an $\mf{sp}_6 \C$- partition can be identified with a vector of the form $$(m_{12}, m^+_{12},m_{13}, m^+_{13}, m_{23}, m^+_{23}, m_1, m_2, m_3 ).$$
Firstly we deal with the case of $\lambda= 2 \omega_3$. We have $c_3=0$, so we set  $m_3=0$ and the Step~1 returns the null vector. For Step~2, we have $c_2=2=c_3+2$ and we are in case $(3)$. We set $m_2=m_3=0$ and $m_{23}=m^+_{23}=1$ obtaining the vector $(0,0,0,0,1,1,0,0,0)$. Finally we have $c_1=c_2+2$ and to perform Step~3 we are again in case~$(3)$, so we set $m_1=m_2=0$, $m_{13}=m^+_{13}=1$ and  $m_{12}=m^+_{12}=1$ and the iteration produces the vector $(1,1,1,1,1,1,0,0,0)$. \\
We want now obtain an $\mf{sp}_6\C$-partitions associated to $2 \rho - 4 \omega_1$. We have that $c_3=c_1=2$ and $c_2=4$. Because $c_3=2$, Step~1 of our construction produces the vector $(0,0,0,0,0,0,0,0,2)$. We have $c_2 = c_3+2$ and then, to perform Step~2, we are in case~$(3)$. We set $m_2=2$ and $m_{23}=m^+_{23}=1$ and we obtain the vector $(0,0,0,0,1,1,0,2,2)$. Finally, because $c_1=c_2-2>0$, at Step~3 we are in case~$(2)$. Observe that $J_2=\{3\}$ and $s=0$, so we set $m_1=m_2=2$ and $m_{12}=m^+_{12}=m_{13}=m^+_{13}=0$, and $(0,0,0,0,1,1,2,2,2)$ is a $\mf{sp}_6\C$-partitions associated to~$2 \rho - 4 \omega_1$.
\end{example}
\begin{example}\label{es:gpartB}
We construct now an admissible $\mf{so}_7\C$-partition associated to the weight $2 \rho - \lambda$ where $\lambda = 4 \omega_1 + 2\omega_3 $. We identify an $\mf{so}_7 \C$- partition $m$ with a vector of the form $$(m_{12}, m^+_{12},m_{13}, m^+_{13}, m_{23}, m^+_{23}, m_1, m_2, m_3 ).$$
In $B_3$ the weight $2 \rho$ has coordinates $(5,3,1)$ with respect to the $\{\varepsilon_i\}$ basis, and then $c_3=0$, $c_2=2$ and $c_1=0$. Consequently we have that $m_{23}=m^+_{23}=1$ are the only non zero coordinates of $m$ and then the algorithm produces the vector $(0,0,0,0,1,1,0,0,0)$.
\end{example}
\textbf{Case B: there exists an even number of odd $c_i$, $c_n$ is even or $c_n$ is odd and $\lambda_n \neq 0$.}
\begin{itemize}
 \item[\emph{Step 1}] Let  $\{\gamma_1 <  \dots < \gamma_{2k}\}$ be the set of indices such that $c_i$ is odd. We pair together the $j$-th and the $k+j$-th index obtaining the set $P=\{(\gamma_1,\, \gamma_{k+1}), \, \dots ,\, (\gamma_k, \, \gamma_{2k})\}$.
 \item[\emph{Step 2}] Construct the weight $\lambda'$ starting from $\lambda$ using the pairs in $P$: if $(\gamma_j, \, \gamma_{j+k})\in P$, set $\lambda'_{\gamma_j}=\lambda_{\gamma_{j}}+1$ and $\lambda'_{\gamma_{j+k}}=\lambda_{\gamma_{j+k}}-1$, otherwise $\lambda'_{\gamma_j}=\lambda_{\gamma_j}$.
\item[\emph{Step 3}] Observe that $\lambda'$ is again a dominant weight smaller than $2 \rho$  and the set $\{c'_i=2|\rho_i|-|\lambda'_i|\}$ is composed only by non negative even integers. Using Case~A, construct an admissible $\g$-partition  $m'=(m'_{ij},\, m'^+_{ij}, \,m'_i)$ associated to $2\rho - \lambda'$.
\item[\emph{Step 4}] If $(\gamma_j, \, \gamma_{j+k})$ is a pair in $P$, we set $m_{\gamma_j \gamma_{j+k}}=m'_{\gamma_j \gamma_{j+k}}+1$, otherwise  $m_{\gamma_j \gamma_{j+k}}=m'_{\gamma_j \gamma_{j+k}}$.
 \end{itemize}
 \begin{rmk}\label{rmk:gpartprop} A $\g$-partition constructed in Case~B has the following properties:
 \begin{enumerate}
 \item  $m_{ij}>1$ only if $(i,j)$ is in $P$;
 \item  $m^+_{i j}$ is different from $0$ only if $m_{i j}\neq 0$. Moreover we have $m_{i j} \leq 2 $ and $m^+_{i j}\leq 1$. In particular $m_{ij} > m^+_{ij}$ if and only if $(i,j)=(\gamma_h,\gamma_{h+k}) \in P$. Analogously, $M(i,j)\neq 0$ if and only if  $i = \gamma_h$ and $j = \gamma_{h+k}$, in that case we have $M(i,j)=1$; 
\item  $m^+_{ij}\neq 0$ only if $m^+_{i+1j} \neq 0$ or if $j=i+1$. Consequently the quantities $R(i,j)=m^+_{i j}-m^+_{i+1 j}$ and $m^+_{i j}-m_{i+1 j}$ are smaller or equal to zero if $j>i+1$;
\item $m_i \neq 0$ only if $m_{i+1} \neq 0$. This implies $m_i-m_{i+1}\leq 0$ for all $i$. Moreover observe that for every $i$ we have  \[m_i=\begin{cases} \leq 1 \mbox{ in type } B,\\ 0 \mbox{ in type } D,\\ \leq 2 \mbox{ in type } C.\end{cases}\]
\item
Because of (2), we have that $S(i,j)=m_{ij}+m^+_{i j}-(m_{i+1 j}+m^+_{i+1 j})$ is always smaller or equal to zero, except if $(i,j)=(\gamma_h, \gamma_{k+h}) \in P$. In this case we have $m_{ij}+m^+_{i j}-(m_{i+1 j}+m^+_{i+1 j})=1$.\end{enumerate} \end{rmk}
\begin{proposition} 
The construction exposed in Case B produces an admissible $\g$-partition associated to~$2 \rho  - \lambda$. \end{proposition}
\proof First of all observe that (3) and (4) in Remark~\ref{rmk:gpartprop}
imply immediately that $ \mathscr{N}_i^{t \; 0}(m)$, $ \mathscr{N}_i^{\overline{t} \; 0}(m)$ and $ \mathscr{N}_i^{n \; 0}(m)$ are all smaller or equal than $1$. 
We want now find an upper bound to $\mathscr{N}_i^{\overline{t} \; 1}(m)$ and $\mathscr{N}_i^{t \; 1}(m)$.
We have to discuss some cases, depending on the parity of~$c_i$ and~$c_{i+1}$. Set 
$P_-:=\{\gamma_1, \dots, \gamma_k\} $ and $ P_+:=\{\gamma_{k+1} \dots \gamma_{2k}\}.$\\
\emph{\textbf{If $c_i$ is even}} By construction in Case A we have that $m_{ij+1}+m^+_{i j+1}=m'_{ij+1}+m'^+_{i j+1}\leq m'_{i+1 j+1}+m'^+_{i+1 j+1}$ for $j\neq i$ and then  $S(i,j+1)=m_{ij+1}+m^+_{i j+1}-(m_{i+1 j+1}+m^+_{i+1 j+1})$ is non positive for every $j> i+1$. Moreover $M(i,j)=0$ for all $j$ and $N(i)\leq 0$. It is immediate to check that $ \mathscr{N}_i^{\overline{t} \; 1}(m) \leq 1$ and $ \mathscr{N}_i^{t \; 1}(m) \leq 1$;\\
\emph{\textbf{If $c_i$ is odd and $i \in P_+$}}, by (5) of Remark~\ref{rmk:gpartprop} we have that $S(i,j)\leq 0$ and $M(i,j)=0$ for every $j$. The  inequalities for $\mathscr{N}_i^{\overline{t} \; 1}(m)$ and $\mathscr{N}_i^{t \; 1}(m)$ are then easily verified;\\
\emph{\textbf{If $c_i$ and $c_{i+1}$ are both odd and $i,i+1 \in P_-$}}, suppose $i=\gamma_h$ (and then $i+1=\gamma_{k+h+1}$). We have $S(i, \gamma_{k+h})=1$ and $S(i, \gamma_{k+h+1})<0$.
It follows that for every $s>i+1$
\begin{equation}\label{brutta}\sum_{j=s}^{n-1}\left[m_{ij+1}+m^+_{i j+1}-(m_{i+1 j+1}+m^+_{i+1 j+1})\right] \leq 0.\end{equation}
An immediate consequence of above inequality and of (3) and (4) of Remark~\ref{rmk:gpartprop} is that $ \mathscr{N}_1^{\overline{t} \; 1}(m)- m^+_{i i+1}\leq 0$. Because of (2) of Remark~\ref{rmk:gpartprop} we have $m^+_{i i+1}\leq 1$ and then $ \mathscr{N}_1^{\overline{t} \; 1}(m) \leq 1$.
Observe now that if $s \geq \gamma_{k+h}$ inequality in (\ref{brutta}) is strict. Moreover $M(i, j) = 1$ only if $j=\gamma_{k+h}$ and we obtain consequently that $ \mathscr{N}_i^{t \; 1}(m) \leq 1$ for every $t$ in Table 1.\\
\emph{\textbf{If $c_i$ and $c_{i+1}$ are both odd, $i \in P_-$ and $i+1\in P_+$}}, observe that $c'_{i} \leq c'_{i+1}$ by Step~2 of the construction in Case~B and this implies that $m^+_{i i+1}=m'^+_{i i+1}=0$. Because $i \in P_-$, we can suppose $i=\gamma_h$ and we recall that $M(i,j)>0$ if and only if $j=\gamma_{k+h}$. Moreover we have \[\sum_{j=s}^{n-1}S(i,j+1) =\begin{cases}\leq 1 & \mbox{ if } s < \gamma_{k+h}\\ \leq 0 & \mbox{ otherwise.}\end{cases}\]
Observe now that $M(i,j)>0$ (in particular it is equal to 1) only if $\sum_{j=s}^{n-1}S(i,j+1) \leq 0$ and the inequalities $\mathscr{N}_i^{\overline{t} \; 1}(m) \leq 1$ and $\mathscr{N}_i^{t \; 1}(m)\leq 1$ are verified; \\
Finally, \emph{\textbf{if $c_i$ is odd, $i=\gamma_h \in P_-$ and $c_{i+1}$ is even}}, we observe again that because of Step~2 of our construction in Case~B, we have $c'_i \leq  c'_{i+1}$ and then  $m_{ii+1}=m^+_{ii+1}=0$. As in the previous case we have \[\sum_{j=s}^{n-1}S(i,j+1) =\begin{cases}\leq 1 & \mbox{ if } s < \gamma_{k+h}\\ \leq 0 & \mbox{ otherwise.}\end{cases}\] 
and $M(i,j)=1$ only if $\sum_{j=s}^{n-1}S(i,j+1) \leq 0$. Check that $\mathscr{N}_i^{\overline{t} \; 1}(m) \leq 1$ and $\mathscr{N}_i^{t \; 1}(m)\leq 1$ in now completely straightforward.\\
It remains to prove that the conditions of Theorem~\ref{BZ} holds for the operators $ \mathscr{L}^s_j(m)$.
Some of these inequalities are trivial by the construction of $m$, in particular $\mathscr{L}_n^{\overline{t}}(m), \mathscr{L}_n^{t}(m) \leq 1$; in fact $m_i \leq 1 $ in type $B$ and $D$, $m_i/2 \leq 1 $ in type $C$ and in our construction we have $M(i,j)\geq 0$ and $ m^+_{ij} \leq 1$ for every $i,j$.
Furthermore observe that $ \mathscr{L}^s_j(m)=\mathscr{L}^{\overline{s-1}}_j(m) - M(s, j)$ and, again because $M(i,j)$ are always non negative,  we reduce to prove that $\mathscr{L}^{\overline{s-1}}_j(m) \leq 1$.
We recall that
 \[\mathscr{L}_j^{\overline{t}}(m)=    \sum_{i=1}^{t } \left( M(i, j+1) - M(i, j)\right)    + m_{t +1\, j+1}. \]
We have four cases:\\
\emph{\textbf{If both $j$ and $j+1$ are not in $P_{+}$}} by (2) of Remark~\ref{rmk:gpartprop} we have $M(i,j)=M(i, j+1)=0$ for all $h$ and for all $j$. Moreover $m_{t+1 j+1}$ is smaller than  $1$  because $j+1 \notin P_{+}$ and the inequality $\mathscr{L}_j^{\overline{t}}(m) \leq 1$ is verified.\\
\emph{\textbf{If $j=\gamma_{k+h} \in P_{+} $ and $j+1 \notin P_{+}$ }}, we have $M(i, j+1)= 0$ for all $i$ and $M(i, j)= 0$ if and only if $i \neq \gamma_h$. Otherwise we have $M(\gamma_h, \gamma_{k+h})=1$. 
This implies that $ \sum_{i=1}^{t } \left( M(i, j+1) - M(i, j) \right)=-1$ if $t\geq \gamma_h$. Otherwise $ \sum_{i=1}^{t } \left( M(i, j+1) - M(i, j) \right) =0$. Moreover we have $m_{t+1 j+1} \leq 1$ because $j+1 \notin P_{+k}$. These conditions immediately imply $\mathscr{L}_j^{\overline{t}}(m) \leq 1$.\\
 \emph{\textbf{If $j \notin P_{+} $ and $j+1=\gamma_{k+h} \in P_{+}$, }} we firstly remark that by construction we have $c'_j < c'_{j+1} + 2$ and then  $m'_{j j+1}=0$ by construction in Case~A. Thus $m_{jj+1}=1$ if $j=\gamma_h$ and zero otherwise. In general,  $m'_{j j+1}=0$ implies  $m'_{i j+1}=0$ for every $i\leq j$ and then $m_{t+1 j+1}=1$ if $t=\gamma_h - 1$ and zero otherwise. Moreover observe that $ M(i,j+1)>0$ (and in particular, it is equal to 1) only if $i=\gamma_h$.
Now we can evaluate the expression $ \sum_{i=1}^{t } \left( M(i,j+1) - M(i, j) \right)$. By our previous observations about the $M(i, j+1)$ and by (2) of Remark~\ref{rmk:gpartprop}, it is equal to $0$ if $t< \gamma_h$ and equal to $1$ if $t \geq \gamma_h$. As observed before, in the last case we have $m_{t+1 j+1}=0$ and it follows easily that $\mathscr{L}_j^{\overline{t}}(m) \leq 1$ holds.\\
 \emph{\textbf{If $j$ and $j+1$ are both in $P_+$, }} we can suppose $j=\gamma_{ k+h}$ and then $j+1=\gamma_{ k + h + 1}$. We consequently have $M(i, j) \neq 0$ if and only if $i=\gamma_h$ and that $M(i, j+1) \neq 0$ if and only if $i=\gamma_{h+1}$. We then obtain that  $ \sum_{i=1}^{t } \left( M(i, j+1) - M(i, j) \right)$ is equal to $-1$ if $\gamma_h \leq t < \gamma_{h+1}$ and $0$ otherwise. If $m_{t+1 j+1}\leq 1$ the inequality $\mathscr{L}_j^{\overline{t}}(m) \leq 1$ is verified. Otherwise, we remark  that $m_{t+1 j+1}=2$ only if $t+1=\gamma_{h+1}$, i.e if $\gamma_h \leq t < \gamma_{h+1}$, but this is exactly the case of $ \sum_{i=1}^{t } \left( M(i, j+1) - M(i, j) \right)=-1$, and again the inequality is checked.\endproof

\begin{example}
In this example we want to construct an admissible $\mf{sp}_8\C$-partition associated to the weight $2 \rho - \lambda$,  where $\lambda= \omega_4$. We recall that $\omega_4$ has coordinates $(1,1,1,1)$ in the $\{\varepsilon_i\}$ basis so we have $c_4=1, c_3=3, c_2=5, c_1=7$. The set of odd indices is $\{1,2,3,4\}$ and $P=\{(1,3), (2,4)\}$. The weight $\lambda'$ is then $(2,2,0,0)$ (i.e. $2 \omega_2$) and by construction in case~A we have that the non zero coordinates of $m'$ are $m'_4=m'_3=m'_2=m'_1=2$, $m'_{34}=m'^+_{34}=1$, $m'_{24}=m'^+_{24}=1$ and 
$m'_{12}=m'^+_{12}=m'_{14}=m'^+_{14}=1$. By our construction in case~B, we have that the $\mf{sp}_8\C$-partitions $m$ associated to the weight $2 \rho - \omega_4$ has the following non zero coordinates: $m_4=m_3=m_2=m_1=2$,  $m_{34}=m^+_{34}=1$, $m_{24}=2, m^+_{24}=1$ and $m_{12}=m^+_{12}=m_{13}=m_{14}=m^+_{14}=1$.
\end{example}
\begin{rmk}
Because of the parity constraint in the dominance order relations in type $C$ and $D$ (c.f.r. Remark~\ref{classmin2rho}), Case A and B cover all the weights appearing in the statement of Theorem~\ref{thm:lexdom} for symplectic and even orthogonal algebras. 
\end{rmk}
Because of previous Remark, in the remaining cases we deal only with algebras of type $B$. In particular, observe that in type $B$ the condition $c_n \neq 0$ is equivalent to assume that $c_n$ is odd and in particular $c_n=1$. Moreover  $c_n \neq 0$ if and only if $\lambda_n=0$.\\
\textbf{Case C: $\lambda_n=0$ and $c_n$ is odd or there exists an odd number of odd $c_i$ and $\lambda_n\neq0$.}  Let $I=\{\gamma_1   \dots < \gamma_{k}\}$ be the set of indices such that $c_i$ is odd.
\begin{itemize}
 \item[\emph{Step 1}] Construct the weight $\lambda'$ setting $\lambda'_i=\lambda_i+1$ if $i \in I$ and $\lambda'_i=\lambda_i$ otherwise. Observe that $\lambda'$ is a dominant weight and it is again smaller than $2 \rho$ and that the set $\{c'_i=2|\rho_i|-|\lambda'_i|\}$ is composed only by non negative even integers.
\item[\emph{Step 2}] Using Case~A, construct an admissible $\g$-partition  $m'=(m'_{ij},\, m'^+_{ij}, \,m'_i)$ associated to $2\rho - \lambda'$. Observe that $\lambda'_n \neq 0$, then $c'_n=0$ and by construction in Case A we have $m'_j=0$ for every $j$.
\item[\emph{Step 3}] Set $m_{ij}=m'_{ij}$ and $m^+_{ij}=m'^+_{ij}$ for every pair of indices $i,j$. Moreover set $m_i=1$ if $i \in I$ and $m_i=m'_i=0$ otherwise. 
 \end{itemize}
\begin{proposition} 
The construction exposed in Case~C produces an admissible $\g$-partition associated to~$2 \rho  - \lambda$. \end{proposition}
\proof Let $m$ be an admissible $\g$-partition associated to~$2 \rho  - \lambda$ constructed using the iterative process exposed in Case~C.  By our construction, $\mathscr{L}_j^{\overline{t}}(m)=\mathscr{L}_j^{\overline{t}}(m')$ and $\mathscr{L}_j^{t}(m)=\mathscr{L}_j^{t}(m')$ for every $t$ and for every $j \neq n$.  Moreover  $\mathscr{N}_i^{t \; 0} (m)=\mathscr{N}_i^{t \; 0} (m')$ for every $t \neq n$. Observe that $m_i \leq 1$ for every $i$ and $M(i,j)=m_{ij}+m^+_{ij}=m'_{ij}+m'^+_{ij}= 0$ for every pair of indices $i,j$ because of construction in Case~A. This implies that $\mathscr{L}_n^{\overline{t}}(m), \mathscr{L}_n^{t}(m) \leq 1$. Observe now that, by Step 3 of construction in case~C we have 
\[\mathscr{N}_i^{n \; 0} (m)-m^+_{ii+1}-m_i+m_{i+1}=\mathscr{N}_i^{n \; 0} (m')-m'^+_{ii+1}-m'_i+m'_{i+1}\]\[\mathscr{N}_i^{t \; 1} (m)-m^+_{ii+1}-m_i+m_{i+1}=\mathscr{N}_i^{t \; 1} (m')-m'^+_{ii+1}-m'_i+m'_{i+1}\] In particular, by construction of $m'$ both expressions $\mathscr{N}_i^{n \; 0} (m')-m'^+_{ii+1}-m'_i+m'_{i+1}$ and $\mathscr{N}_i^{t \; 1} (m')-m'^+_{ii+1}-m'_i+m'_{i+1}$ are smaller or equal than 0. To prove that $\mathscr{N}_i^{n \; 0} (m), \mathscr{N}_i^{t \; 1} (m) \leq 1$ it is enough to show that $m^+_{ii+1}+m_i-m_{i+1} \leq 1$ for every $i$.
We remark that in our construction $m^+_{ii+1}$ is always smaller than 1 and $m_i\neq 0$ if and only if $i \in I$. Now, if $c_i$ is even, the inequality $m^+_{ii+1}+m_i-m_{i+1} \leq 1$ comes directly from the fact that $m_i\leq m_{i+1}$. If $i$ and $i+1$ are both odd, then $m_i=m_{i+1}=1$ and 
 $m^+_{ii+1}+m_i-m_{i+1} \leq 1$ is satisfied. Finally, if $c_i$ is odd and $c_{i+1}$ is even, observe that $c_{i} \leq c_{i+1} + 1$  by parity constraint and then $c'_i \leq c'_{i+1}$ by Step~1 of construction in case~C. This implies, by construction of $m'$ and by Step~3 in case~C, that $m^+_{ii+1}=m'^+_{ii+1}=0$ and again we obtained  $m^+_{ii+1}+m_i-m_{i+1} \leq 1$.\endproof
\begin{rmk}\label{subsec:comments}
In type $B_n$, the construction of admissible $\g$ partitions exposed in Case~C works also in Case~B.
We privileged the procedure exposed in Case~B to underline the uniform construction in all the classical cases. \end{rmk}
\begin{example}
In this example we construct admissible $\mf{so}_7\C$-partition associated to the weight $2 \rho - \lambda$ where $\lambda = 4 \omega_1 $. 
We have $c_3=1$, $c_2=3$ and $c_1=1$. The weight $\lambda'$ constructed as in Step~1 of case~C has coordinates $(5,1,1)$ (i.e. $\lambda'=4\omega_1 + 2\omega_3$). We have just constructed in Example~\ref{es:gpartB} an $\mf{so}_7\C$-partition  $m'$ associated to $2 \rho - \lambda'$. In particular we obtained $m'=(0,0,0,0,1,1,0,0,0)$. By Step~3 of construction in case~C, we then obtain that $m=(0,0,0,0,1,1,1,1,1)$ is an $\mf{so}_7\C$-partition associated to~$2 \rho - 4 \omega_1$.
\end{example}
\subsection{A Conjecture about Exterior Algebra $\Lambda V_{\theta_s}$}
If $\g$ is not simply laced, we propose here an analogous of Kostant Conjecture about irreducible representations appearing in the exterior algebra over the little adjoint representation. We are motivated by two recent works that highlight some interesting aspects of the structure of 
$\Lambda V_{\theta_s}$ as $\g$-representation. The first one is an article of I. Ademehin~\cite{IA}, dealing with the graded multiplicities of trivial and little adjoint representation in $\Lambda V_{\theta_s}$. The results contained in~\cite{IA} are in some sense very similar to the classical ones about exterior algebra over $\g$ and we think that a further investigation about the structure of $\Lambda V_{\theta_s}$ could lead to some very interesting results.
Our second motivating paper is an article of Panyushev \cite{Pan}, where the following theorem is proved in the more generic context of orthogonal isotropy representations.
 \begin{theorem}[Panyushev \cite{Pan}, Theorem 2.9]\label{Pan} Let $\g$ be a non simply laced algebra of type $B$, $C$ and $F_4$. Let $\theta_s$ be the short dominant root of $\g$, then
  \[\Lambda V_{\theta_s} \simeq 2^{|\Delta_s|}\left(V_{\rho_s} \otimes V_{\rho_s}\right)
   \]
where $\Delta_s$ is the set of short simple roots and $\rho_s$ is half the sum of positive short roots.
 \end{theorem}
Analogously to the case exterior algebra over adjoint representation, we formulate the following conjecture:
\begin{conj}\label{rhos}Let $\g$ be a non simply laced simple Lie algebra.
 $V_\lambda$ is an irreducible component of $\Lambda V_{\theta_s}$ if and only if $\lambda \leq 2 \rho_s$.
\end{conj}
By Theorem \ref{Pan}, Conjecture \ref{rhos} can be restated as
\begin{conj}\label{Conjrhoshort}Let $\g$ be a non simply laced simple Lie algebra.
 $V_\lambda$ is an irreducible component of $V_{\rho_s} \otimes V_{\rho_s}$ if and only if $\lambda \leq 2 \rho_s$.
\end{conj}
The Conjecture \ref{Conjrhoshort} can be easily proved for case $B_n$ using elementary representation theory.
We checked the conjecture also using Berenstein and Zelevinsky polytope associated to  $V_{\rho_s}\otimes V_{\rho_s}$.
Moreover, we proved it for exceptional cases $F_4$ and $G_2$ by direct computations. The conjecture remains open only in type $C$, where it seems that combinatorics of short roots and weights is linked to the Kostant conjecture in type $D$.

\section{Generalized Exponents and Macdonald Kernels}\label{sec:exp}
We give here an overview of theory of generalized exponents for representations of Lie algebras, following the results exposed in \cite{KostPoly}.
\begin{theorem}[Kostant~\cite{KostPoly}, Theorem 0.11]\label{GradedS}
 The module $\mathrm{Hom}_\g \left(V_\lambda, S(\g)\right)$ is a free $S(\g)^\g$-module of dimension $\dim V_\lambda^0$. 
\end{theorem}
Let $n$ be the dimension of $V_\lambda^0$ and let $f_1, \dots, f_n$ be any set of homogeneous generators of $\mathrm{Hom}_\g \left(V_\lambda, S(\g)\right)$  as $S(\g)^\g$-module. Up to relabeling the polynomials $f_i$, it is possible to suppose that~$\mathrm{deg}f_i \leq \mathrm{deg}f_{i+1}$ for every $i$. Set $m_i(\lambda)=\mathrm{deg}f_i$.
\begin{defn}[c.f.r \cite{Broer}, Introduction]
 The integers $m_1(\lambda),  \dots , m_n(\lambda)$ are the \emph{generalized exponents} of the representation $V_\lambda$.
\end{defn}
Generalized exponents have also an interpretation in therms of $W$-representation on the zero weight space $V_\lambda^0$.
Let $c \in W$ be a Coxeter - Killing transformation, i.e. $c = s_{\alpha_1} \dots s_{\alpha_n}$ where $s_{\alpha_i}$ is the simple reflection associated to the $i$-th simple root.
The action of $\g$ on $V_\lambda$ induces a representation $\rho_\lambda: W \rightarrow \mathrm{End}(V_\lambda^0)$. The element $\rho_\lambda(c)$ acts diagonally on $V_\lambda^0$ with eigenvalues $\gamma_j = \mathrm{exp}{\frac{2i \pi m_j(\lambda) }{h}} $, where $h$ is the Coxeter number. 
\begin{example}
Consider $\g$ acting on itself by the adjoint action. Such an action induces the reflection representation of $W$ on $\h$. The generalized exponents for this representation coincides with the classical exponents of $\g$.
\end{example}
\begin{example}
Let $\g$ be a not simply laced simple Lie algebra and consider its little adjoint representation $V_{\theta_s}$.The generalized exponents associated to $V_{\theta_s}$ are the \emph{short exponents} of $\g$.
\end{example}
Consider now the generating polynomial of generalized exponents defined by the formula
  \[E_\lambda(t)= \sum_{i=1}^{\dim V_\lambda^0} t^{m_i(\lambda)}.\]
Theorem \ref{GradedS} translates naturally into the following remarkable factorization of generating series of graded multiplicities:
\begin{equation*} P(V_\lambda, S(\g),  t)= E_\lambda(t)\prod_{i=1}^n (1-t^{e_i+1})^{-1}.\end{equation*}
Determining the graded multiplicities in the symmetric algebra is then deeply linked to determining the generalized exponents of $V_\lambda$. In particular the problem of finding explicit formulae for the polynomials $E_\lambda(t)$ turns out to be very interesting both from a combinatorial and from a representation theoretic point of view. For Lie algebras of type $A$ a combinatorial description of generalized exponents is given in \cite{LaSchu}, \cite{LaLT}. For other classical algebras, the combinatorics of generalized exponents is less explicit, and closed formulae are available only in special cases (see \cite{GUP}, \cite{Ion1}, \cite{Ion2}, \cite{Ion3}, \cite{LL}, \cite{Sh}).
\begin{rmk}\label{rem:explge}
Values of classical exponents, and of short exponents in the case of non simply laced algebras, are well known (see \cite{Stembridge}, Table 4.1). Explicit formulae for $E_\theta(t)$ and $E_{\theta_s}(t)$ can be consequently computed in the classical cases:
\begin{equation}
   E_\theta(t)=  \left(n+1\right)_t\qquad \mbox{ Type } A_n
\end{equation}
\begin{equation}
 E_\theta(t)=t\left(n\right)_{t^2} \qquad E_{\theta_s}(t)=t^n\qquad \mbox{ Type } B_n
\end{equation}
\begin{equation}
     E_\theta(t)= t\left(n\right)_{t^2} \qquad E_{\theta_s}(t)=t^2\left(n-1\right)_{t^2} \qquad \mbox{ Type } C_n
\end{equation}
\begin{equation}
     E_\theta(t)=\left(n\right)_{t^2} \frac{t(t^{n-2}+1)}{(t^n+1)} \qquad \mbox{ Type } D_n
\end{equation}
\end{rmk}
\subsection{Macdonald Kernels} We recall now some tools, introduced by Stembridge in \cite{Stembridge}, useful to produce effective computations.
\begin{defn}[c.f.r. \cite{Stembridge}, Section 1.1]
Let $\Z\langle \Pi \rangle:=\Z\{e^\lambda,\,\lambda \in \Pi\}$ denote the group ring generated by $\Pi$. The Macdonald Kernel of $\g$ is the formal series  $\Delta(q,t) \in \Z\langle \Pi \rangle\left[[q,t]\right]$ defined by the formula
  \[\Delta(q,t):= \prod_{i\geq 0 } \left( \frac{1-q^{i+1}}{1-tq^i}\right)^{\mathrm{rk}\g} \cdot \prod_{i \geq 0} \prod_{\alpha \in \Phi} \frac{1-q^{i+1}e^\alpha}{1-tq^ie^\alpha}.\]
\end{defn}
The Macdonald kernels specialize to the graded character of exterior algebra of adjoint representation, when evaluated at $(q,t)=(-q,q^2)$, and to the graded character of the symmetric algebra over $\g$ when evaluated at $(q,t)=(0,t)$ (c.f.r.~\cite{Stembridge}, Section 1.2). Observe now that $\Delta(q,t)$ is $W$-invariant; this implies that it can be expanded in terms of characters of irreducible representations, obtaining an expression of the form
\[\Delta(q,t)=\sum_{\mu\in  \Pi^+}  C_\mu (q,t) \chi(\mu),\]
for certain formal series $C_\mu (q,t)$, indexed by dominant weights of $\g$. In particular, when specialized at $(q,t)=(-q,q^2)$ and at $(q,t)=(0,t)$,  the formal series $C_\mu (q,t)$ gives the graded multiplicities of the representation $V_\mu$ in the exterior algebra and in the symmetric algebra respectively.
\begin{rmk}\label{rmk:genratio}By Theorem \ref{GradedS}, the polynomial $E_\lambda(t)$ can be computed by determining the ratio $C_\lambda (q,t)/C_0(q,t)$ and evaluating it at $(0,t)$.\end{rmk}
In  \cite{Stembridge} Stembridge proves that the formal series $C_\mu(q,t)$ satisfy some recurrences, reducing the problem of their explicit computation to solving a linear system of equations with coefficients in 
$\C[q^{\pm 1},t^{\pm 1}]$.\\
We recall that it is possible to extend the definition of $C_\mu (q,t)$ to any weight $\mu$ by setting
 \begin{equation*}\label{RedClambda}
C_\mu(q,t) =
\left\{
	\begin{array}{lll}
	0 & \mbox{if } \mu + \rho \mbox{ is not regular}, \\
    (-1)^{l(\sigma)}C_\lambda(q,t)  & \mbox{if } \sigma(\mu + \rho) = \lambda+\rho \; , \lambda \in \Pi^+, \sigma \in W .\\
	\end{array}
\right.
\end{equation*}
For short, if there exists $\sigma$ such that $\sigma(\mu + \rho)=\lambda + \rho$, with $\lambda \in \Pi^+$, we say that the weight $\mu$ is \emph{conjugated to $\lambda$ by $\sigma$}. Moreover, if $\mu$ is conjugated to $\lambda$ by $\sigma$, we say that $(-1)^{l(\sigma)}C_\lambda(q,t)$ is the \emph{reduced form} of $C_\mu(q,t)$. Sometimes a precise information about the sign of $\sigma$ is not needed in our reasoning; in this case we shortly say that \emph{$\mu$ is conjugated to $\lambda$}.
\begin{theorem}[ Minuscule Recurrence, \cite{Stembridge}, Formula (5.14)]
Fix a dominant weight $\lambda$ and let $\omega $ be a minuscule coweight (i.e. $(\omega, \alpha ) \in \{0, \pm1\} $ for every positive root $\alpha$), then the following relation holds:
 \begin{equation}\label{ricorsionemin}
  \sum_{i=1}^k C_{w_i \lambda}(q,t) \left(\sum_{\psi \in O_\omega} \left(t^{-(\rho, w_i \psi)}-q^{(\lambda, \omega)}t^{(\rho, w_i\psi)}\right)\right)=0.
 \end{equation}
where, denoting by $W_\lambda$ the stabilizer of $\lambda$ in $W$, the $w_1, \dots , w_k$ are minimal coset representatives of $W/W_\lambda$ and $O_\omega$ is the orbit $W_\lambda \cdot \omega$.\end{theorem}
\begin{rmk}\label{rmk:reducedsmaller}
Observe that if $\lambda$ and $\mu$ are dominant weights and $w \lambda$ is conjugated to $\mu$, then $\mu < \lambda$. As a consequence, if we write the $ C_{w_i \lambda}(q,t)$ appearing in Formula~(\ref{ricorsionemin}) in their reduced form, in the minuscule recurrence there appear only  $C_{\mu}(q,t)$ with $\mu$ dominant and smaller than~$\lambda$.
\end{rmk}
\section{Small Representations}\label{sec:small}
The aim of this Section is to present closed formulae for generalized exponents of certain small representations in type $B$, $C$ and $D$. In Table~\ref{tab:2} we list the weights of non trivial small representations in these three cases. More precisely, in Theorem~\ref{thm:geC}, Theorem~\ref{thm:geB} and Theorem~\ref{thm:geD} we provide closed expressions for the polynomials of generalized exponents for small representations that are indexed by fundamental weights.
 \begin{table}[ht!]
  \begin{center}
   \caption{Weights of small representation for type $B, C$ and $D$.}\label{tab:2}
   \begin{tabular}{|c|c|c|c|}
   \hline
    \textbf{Type B} & \textbf{Type C} & \textbf{Type D, $n$ even} & \textbf{Type D, $n$ odd} \\
    \hline
     $\omega_i$, $i<n$& $\omega_{2k}$&$\omega_{2k}$&$\omega_{2k}$\\
     $2 \omega_n$ & $\omega_1+\omega_{2k+1} $&$2\omega_{n-1}, \, 2\omega_{n}$&$\omega_{n-1}+\omega_{n}$ \\
  & & $\omega_1 + \omega_{2i+1}$& $\omega_1 + \omega_{2i+1}$\\
   & & $\omega_1 + \omega_{n-1}+\omega_{n}$&$\omega_1 + 2\omega_{n-1}, \, \omega_1 + 2\omega_{n-1} $\\
\hline
   \end{tabular}
\end{center}
\end{table}
 \begin{theorem}\label{thm:geC}Let $\lambda$ be a small weight of the form $\lambda =\omega_{2k}$ for the simple Lie algebra of type~$C_n$. Then:
 \begin{equation}\label{eqn:geC}
E_{\lambda}(t)=\frac{t^{2k}(n-2k+1)_{t^2}}{(n-k+1)_{t^2}}\binom{n}{k}_{t^2}.
 \end{equation}
\end{theorem}
\begin{theorem}\label{thm:geB}The polynomials of generalized exponents for small weight for the simple Lie algebra of type $B_n$ have the following closed expressions:
\begin{equation}\label{eqn:geB}E_{\omega_{2k}}(t)=t^k\binom{n}{k}_{t^2}, \quad E_{\omega_{2k+1}}(t)=t^{n-k}\binom{n}{k}_{t^2},\quad E_{2\omega_{n}}(t)=t^{n-\lfloor\frac{n}{2}\rfloor}\binom{n}{\lfloor\frac{n}{2}\rfloor}_{t^2}, \end{equation}
\end{theorem}
\begin{theorem}\label{thm:geD} Let $\lambda$ be a small weight of the form $\lambda=\omega_{2k}$ for the simple Lie algebra of type~$D_n$, then:
\begin{equation}\label{eqn:geD1}E_{\omega_{2k}}(t)=t^k\frac{(t^{n-2k}+1)}{(t^n+1)}\binom{n}{k}_{t^2},\end{equation}
Moreover, if $n$ is odd, the weight $\omega_{n-1}+\omega_n$ is small and we have:
\begin{equation}\label{eqn:geD2} E_{\omega_{n-1}+\omega_n}(t)=\frac{t^{\lfloor\frac{n}{2}\rfloor}(t+1)}{(t^n+1)}\binom{n}{\lfloor\frac{n}{2} \rfloor}_{t^2}, \end{equation}
Finally, if $n$ is even, $2\omega_{n-1}$ and $2\omega_n$ are small and the following formulae hold:\begin{equation}\label{eqn:geD3} E_{2\omega_{n-1}}(t)=E_{2\omega_{n}}(t)=\frac{t^{\frac{n}{2}}}{(t^n+1)}\binom{n}{\frac{n}{2}}_{t^2},\end{equation}
\end{theorem}
The proof of above theorems use and iterative reasoning, based on the fact that if $\lambda$ is small and $\lambda' < \lambda$, then $\lambda'$ is small. In particular, a minimal non zero small weight is a dominant root. 
The following remark provides the base step for our computations.
\begin{rmk}\label{rmk:base}
The following formulae are proved in \cite{Stembridge}, Theorem 4.1:
 \begin{equation}\label{eqn:ad}
 C_{\theta}(q,t)=\frac{t-q}{t-qt^{h}}E_\theta(t) C_0(q,t), \qquad  C_{\theta_s}(q,t)=\frac{t-q}{t-qt^{h}}E_{\theta_s}(t) C_0(q,t)
  \end{equation}
\end{rmk}
\subsection{Proof of Formulae in Type C} In this section we give a proof of Theorem~\ref{thm:geC}. Observe that such a theorem hold true if $k=1$ because of Remark~\ref{rem:explge}.
 In~\cite{SDT}, Theorem 5.5, the following iterative formula for $C_\lambda (q,t)$ is achieved for weights of the form $\lambda=\omega_{2k}$:
\begin{equation}\label{odiolatex}
          C_{\omega_{2(k+1)}}(q,t)= \frac{(t^{2(n-2k-1)}-1)(t^{2(n-k+1)}-1)(1-qt^{2k-1})t^2}{(t^{2(n-2k+1)}-1)(t^{2(k+1)}-1)(1-qt^{2(n-k)-1})} C_{\omega_{2k}}(q,t).
            \end{equation}
By Remark~\ref{rmk:genratio}, evaluating Equation~(\ref{odiolatex}) at $q=0$ we obtain a recursive relation between $E_{\omega_{2(k+1)}}(t)$ and $E_{\omega_{2k}}(t)$. Using Remark~\ref{rmk:base} for the base step, we obtain by induction that
\begin{align*}
  E_{\omega_{2(k+1)}}(t)&=\frac{t^2(t^{2(n-2k-1)}-1)(t^{2(n-k+1)}-1)}{(t^{2(n-2k+1)}-1)(t^{2(k+1)}-1)}E_{\omega_{2(k)}}(t)\\
   &=\frac{t^2(t^{2(n-2k-1)}-1)(t^{2(n-k+1)}-1)}{(t^{2(n-2k+1)}-1)(t^{2(k+1)}-1)}\frac{t^{2k}(n-2k+1)_{t^2}}{(n-k+1)_{t^2}}\binom{n}{k}_{t^2}\\
   &=\frac{t^{2(k+1)}(n-2k-1)_{t^2}}{(n-k)_{t^2}}\binom{n}{k+1}_{t^2}
 \end{align*}
and Equation~(\ref{eqn:geC}) is proved. We remark that similar but more complicated formulae can be obtained for the other small weights in type $C$ by making explicit the coefficients of the equations in \cite{SDT}, Section 5.3.
\subsection{Proof of Formulae in Type B}\label{subsec:ReederB}
 In type $B_n$ the unique minuscule coweight is $\varepsilon_1$. In Formula~(\ref{ricorsionemin}) we choose $\omega=\varepsilon_1$ and $\lambda=\omega_k$ with $k<n$ or $\lambda=2 \omega_n$. The stabilizer $W_{\omega_i}$ is isomorphic to $S_{i}\times B_{n-i}$ and $W_{\omega_i}(\varepsilon_1)=\{\varepsilon_1,\dots, \varepsilon_i\}$.
Analogously, $W_{2\omega_n}$ is isomorphic to $S_{n}$ and $W_{2\omega_n}(\varepsilon_1)=\{\varepsilon_1,\dots, \varepsilon_n\}$. Writing all the $C_\mu(q,t)$ in their reduced form the recurrence can be rewritten as
\begin{equation}\label{redSt2}
\sum_{\mu \leq \lambda}\Lambda_\mu^{\lambda,n}(q,t)C_\mu(q,t) = 0,
\end{equation}
for certain coefficients $\Lambda_\mu^{\lambda,n}(q,t)$. We will refer to this form of the Formula~(\ref{ricorsionemin}) as the \emph{reduced recurrence for $C_\lambda(q,t)$}. We recall that if $\lambda$ is a small, then a dominant weight $\mu$ smaller than $\lambda$ in the dominance order is again small.  In particular the weights smaller than $\omega_k$ are of the form $\omega_h$ with $h<k$. From now, we denote by $C_h$ the formal series $C_{\omega_h}(q,t)$ and by $\Lambda_h^{k,n}$ its coefficient $\Lambda_\mu^{\lambda,n}(q,t)$ in the recurrence for $\lambda=\omega_k$. Moreover, if $\lambda=2 \omega_n$, we use the notation $C_n$ for $C_{2\omega_n}(q,t)$ and the coefficient of $C_h(q,t)$ in the relative recurrence will be $ \Lambda_h^{n,n}$. Reasoning as in \cite{SDT} and aiming to simplify recurrence (\ref{redSt2}), we want now expand recursively  the coefficients $\Lambda_h^{k,n}$.
\begin{rmk}\label{rmk:contraction}Using explicit realization of fundamental weights and of $\rho$ as presented in Section~\ref{subsec:weights}, it is possible to check  that a weight of the form $w \omega_k$ in $B_n$, with $w \in W$, is conjugated to $\omega_h$ only if $w \omega_k = \varepsilon_1 + \dots + \varepsilon_h + \nu$, where $\nu$ has the first $h$ coordinates equal to 0 when written in the $\{\varepsilon_i\}$ basis. 
In our realization, the root system $B_{n-h}$ can be identified as the root subsystem of $B_n$ given by vectors of the form $\{\pm \varepsilon_i\pm \varepsilon_j\}_{i<j\leq n-h} \cup \{\pm \varepsilon_1, \, \dots \, , \pm \varepsilon_{n-h}\}$ . This identification corresponds to the immersion of $B_{n-h}$ into $B_n$ induced by the immersion of the associated Dynkin diagrams. Under this identification $\nu$ can be thought of as weight of the form $w' \omega_{k-h}$ conjugated to 0 in $B_{n-h}$.
\end{rmk}
\begin{example}
Consider the weight $\omega_4$ in $B_6$ and let $\epsilon_3 $ the element of $ W$ that acts as the sign change on the $3$-rd coordinate. Then $\epsilon_3 \omega_4 = \omega_2 + \nu$ there $\nu=-\varepsilon_3 + \varepsilon_4$. Observe that $\nu + \rho$ has coordinates $(\frac{11}{2},\frac{9}{2},\frac{5}{2},\frac{7}{2},\frac{3}{2},\frac{1}{2})$ in terms of the $\{ \varepsilon_i \}$ basis. In particular $s_3 (\nu+\rho) = \rho$ and $\nu$ is conjugated to zero.
\end{example}
\begin{rmk}\label{rem:Lambdak}
It is immediate to check that $w \lambda$ is conjugated to $\lambda$ if and only if $w \in W_\lambda$. In particular this implies that, if $\lambda = \varepsilon_1 + \dots + \varepsilon_k$, then
\[\Lambda_k^{k,n}=\sum_{i=1}^k \frac{1-qt^{2n-2j+1}}{t^{\frac{2n-2j+1}{2}}}=\frac{(1-qt^{2n-k})(t^k-1)}{t^{\frac{2n-1}{2}}(t-1)}\]
\end{rmk}
The following Lemma, proved in \cite{SDT}, characterize the set  $\Omega^{k,n}_0$ of weights of the form $w \omega_k$ conjugated to 0 in $B_n$. In particular, it is possible to describe explicitly the coordinates of a weight in $\Omega^{k,n}_0$  with respect to the basis $\{\varepsilon_1, \dots, \varepsilon_n\}$.
 \begin{lemma}[\cite{SDT}, Lemma 4.1]\label{zeroconj} Let $w \in W$ be such that $w \omega_k \in \Omega^{k,n}_0$, then:
 \begin{itemize}
  \item \emph{if $k$ is even},  $w \omega_k$ has all the coordinates equal to zero except for $k/2$ pairs of consecutive coordinates of the form $(-1,1)$ and it is conjugated to $0$ by $\sigma \in W$ of length~${k/2}$.
  \item  \emph{if $k$ is odd} then $w \omega_k$  has all the coordinates equal to zero, except for a choice of $(k-1)/2$ pairs of coordinates equal to $(-1,1)$ and for the last one that must be equal to $-1$. In this case $w \omega_k$ is conjugated to $0$ by by $\sigma \in W$ of length~${(k-1)/2+1}$.
 \end{itemize}
 \end{lemma}
The cardinality of $\Omega^{k,n}_0$ can be explicitly computed as a consequence of Lemma~\ref{zeroconj}:
  \begin{equation*} 
|\Omega^{k,n}_0| =
	\begin{cases}
	 \binom{n-\frac{k}{2}}{\frac{k}{2}} & \mbox{if } k \mbox{ is even}, \\
 \binom{n-\frac{k-1}{2}-1}{\frac{k-1}{2}}&  \mbox{if } k  \mbox{ is odd.}
	\end{cases}
	\end{equation*}
	Set
\[p(n,q,t)= t^{\frac{2n-1}{2}}-qt^{-\frac{2n-1}{2}}+t^{-\frac{2n-3}{2}}-qt^{{\frac{2n-3}{2}}}=\frac{(t-q)(1+t^{2n-2})}{t^{\frac{2n-1}{2}}}.\]
\begin{lemma}\label{lem:expansion} The following relations between the coefficients $\Lambda^{k,n}_h$ hold for $h<k$:
  \begin{equation}\label{eqn:expansion1}\Lambda^{k,n}_h= (-1)^s\Lambda^{h,n}_h \binom{n-k+s}{s} + \Lambda^{k-h,n-h}_0 \qquad  \mbox{ if } \; k-h=2s,\end{equation}
   \begin{equation}\label{eqn:expansion1bis}\Lambda^{k,n}_h= (-1)^{s+1}\Lambda^{h,n}_h \binom{n-k+s}{s} + \Lambda^{k-h,n-h}_0 \qquad  \mbox{ if } \; k-h=2s+1\end{equation}
   \begin{equation}\label{eqn:expansion2}\Lambda^{2s,n}_0= (-1)^s p(n,q,t) \binom{n-s-1}{s-1} - \Lambda^{2s-2,n-2}_0 + \Lambda^{2s, n-1}_0, \end{equation}
  \begin{equation}\label{eqn:expansion2bis}\Lambda^{2s+1,n}_0= (-1)^{s+1} p(n,q,t) \binom{n-s-2}{s-1} - \Lambda^{2s-1,n-2}_0 + \Lambda^{2s+1, n-1}_0,\end{equation}
\end{lemma}
\proof 
Equations (\ref{eqn:expansion1}) and (\ref{eqn:expansion1bis}) are direct consequences of Remark~\ref{rmk:contraction}, where it is observed that a weight $w \omega_k$ gives a contribution to $\Lambda_h^{k,n}$ if it is of the form $ \varepsilon_1 + \dots + \varepsilon_h + \nu$, where $\nu$ can be thought as weight $w' \omega_{k-h}$ conjugated to 0 in $B_{n-h}$. 
Moreover, Equations (\ref{eqn:expansion2}) and (\ref{eqn:expansion2bis}) can be obtained observing that a weight $\nu =w \omega_k$ contributing to $\Lambda_0^{k,n}$ is of the form  $-\varepsilon_1+\varepsilon_2+\nu'$ with $\nu'$ conjugated to $0$ in $B_{n-2}$
or of the form $\nu=(0, \nu_2, \dots, \nu_n)$ in $\{\varepsilon_i\}$-expansion, where $\nu'=(\nu_2, \dots, \nu_n)$ is a weight conjugated to $0$ in $B_{n-1}$. 
\endproof
The above relations enable to simplify considerably the computations needed to prove our formulae. In particular they are crucial for the proof, contained in Section~\ref{subsec:fTh}), of the following theorem.
We denote by $R_i$ and by $R_n$ the reduced recurrences for $C_\lambda(q,t)$, with $\lambda=\omega_i$ and $\lambda=2 \omega_n$ respectively. 
\begin{theorem}\label{thm:Bsimply}There exists a family of integers $\{A^{k,n}_i\}_{i\leq k}$ such that 
\begin{equation}\label{eqn:general}
  \sum_{i=1}^k A^{k,n}_i R_i = \Lambda^{k,n}_k C_k + \Gamma_0^{1,n-k+1}C_{k-1} + \sum_{i=1}^{} \Gamma_0^{2, n-k+i+1}C_{k-2i} + \sum_{i=2}\Gamma_0^{2,i}C_{k-2i+1} = 0
\end{equation}
where the coefficients $\Lambda^{k,n}_k$ and $\Gamma^{i,n}_0$ are defined by the formulae

\[\Gamma_0^{1,n}=\Lambda_0^{1,n}=-\frac{(t-q)t^{n-1}}{t^{\frac{2n-1}{2}}}\qquad
\Gamma^{2,n}_0=\Gamma^{2,n-1}_0-p(n,q,t) = -\frac{(t-q)(t^{2n-1}-1)}{t^{\frac{2n-1}{2}}(t-1)}\]
\end{theorem}
Observe that specializing the Equation (\ref{eqn:general}) at $(q,t)\rightarrow (-q, q^2)$ one obtains the equation of \cite{SDT}, Proposition 4.3  used to prove Reeder's Conjecture in type $B$. 
\begin{rmk}\label{rem:ratioexp}Dividing Equation (\ref{eqn:general}) by $C_0(q,t)$ we obtain a recursive relation between the formal series
$\overline{C}_\mu(q,t)=C_\mu(q,t)/C_0(q,t)$. We recall that $\overline{C}_\mu(0,t)=E_\mu(t)$ by Remark~\ref{rmk:genratio}, and consequently the specialization at $q \rightarrow 0$ of Equation~(\ref{eqn:general}) leads to a recursive relation between polynomials of generalized exponents of small representations. \end{rmk}
\proof (of Theorem~\ref{thm:geB}) We denote by $E_h$ and $E_n$ the polynomials $E_{\omega_h}(t)$ and $E_{2\omega_n}(t)$ respectively.
Set
\[b_i=-t^{n-i+1}(t^{2i-1}-1), \qquad c_k= t^k-1\]
Because of Remark~\ref{rem:ratioexp}, evaluating Equation~(\ref{eqn:general}) at $q=0$ and multiplying it by $t^{\frac{2n-1}{2}}(t-1)$ it is possible to obtain the relation
\begin{equation*}
c_k E_k + \sum_{i=1}^{\lfloor\frac{k}{2}\rfloor} b_{n-k+i+1}E_{k-2i} + \sum_{i=1}^{\lfloor\frac{k+1}{2}\rfloor} b_i E_{k-2i+1} = 0.
\end{equation*}
We want now prove our formulae by induction. The base step comes from Remark~\ref{rmk:base} and by formulae contained in Remark~\ref{rem:explge}. 
Consequently, for the inductive step we have to prove that the two identities
\begin{align*}
(t^{2s}-1) t^{s} \binom{n}{s}_{t^2}  &= (t^{2s}-1) E_{2s}\\&= -\sum_{j=0}^{s-1} b_{n-s-j+1}E_{2j} - \sum_{j=0}^{s-1} b_{s-j} E_{2j+1}\\ &= \sum_{j=0}^{s-1} t^{s+j}(t^{2(n-s-j)+1}-1)E_{2j} + \sum_{j=0}^{s-1} t^{n-s+j+1}(t^{2(s-j)+1}-1) E_{2j+1}\\
&=\sum_{j=0}^{s-1} t^{s+j}(t^{2(n-s-j)+1}-1)t^j \binom{n}{j}_{t^2} + \sum_{j=0}^{s-1} t^{n-s+j+1}(t^{2(s-j)+1}-1) t^{n-j} \binom{n}{j}_{t^2}\\&=\sum_{j=0}^{s-1} \left[t^{s+2j}(t^{2(n-s-j)+1}-1) + t^{2n-s+1}(t^{2(s-j)+1}-1) \right] \binom{n}{j}_{t^2}\\&=t^s\sum_{j=0}^{s-1} t^{2j}\left(t^{2(n-2j)}-1\right) \binom{n}{j}_{t^2},
\end{align*}
\begin{align*}
(t^{2s+1}-1) t^{n-s} \binom{n}{s}_{t^2} +b_1 E_{2s}&=(t^{2s+1}-1) E_{2s+1} +b_1 E_{2s}\\&= -\sum_{j=0}^{s-1} b_{s-j+1}E_{2j} - \sum_{j=0}^{s-1} b_{n-s-j} E_{2j+1}\\
&=\sum_{j=0}^{s-1} \left[t^{n+s+1}(t^{2(n-s-j)+1}-1) + t^{n-s+2j}(t^{2(s-j)+1}-1) \right] \binom{n}{j}_{t^2}\\&=t^{n-s}\sum_{j=0}^{s-1} t^{2j}\left(t^{2(n-2j)}-1\right) \binom{n}{j}_{t^2}.
\end{align*}
hold.
Observe now that
\begin{align*}
(t^{2s+1}-1) t^{n-s} \binom{n}{s}_{t^2} +b_1 E_{2s}&= (t^{2s+1}-1) t^{n-s} \binom{n}{s}_{t^2}+b_1 t^{s} \binom{n}{s}_{t^2} \\
&=\left[t^{n-s}(t^{2s+1}-1)-t^{n+s}(t-1)\right] \binom{n}{s}_{t^2}\\
&=t^{n-s}(t^{2s}-1)\binom{n}{s}_{t^2},
\end{align*}
and proving the two identities reduces in both cases to prove that for $s>0$ we have
\begin{equation}\label{eqn:binomials}
(t^{s}-1)\binom{n}{s}_{t}=\sum_{j=0}^{s-1} t^{j}\left(t^{n-2j}-1\right) \binom{n}{j}_{t}.
\end{equation}
Identity~(\ref{eqn:binomials}) can be easily shown by induction on $s$. The case $s=1$ is trivial and 
\begin{align*}
 (t^{s+1}-1)\binom{n}{s+1}_{t}&=\sum_{j=0}^{s} t^{j}\left(t^{n-2j}-1\right) \binom{n}{j}_{t}\\&= t^{s}(t^{n-2s}-1)\binom{n}{s}_t + \sum_{j=0}^{s-1} t^{j}\left(t^{n-2j}-1\right) \binom{n}{j}_{t}\\&=(t^{n-s}-1)\binom{n}{s}_t,
\end{align*}
The equality now holds because
\[\binom{n}{s+1}_{t}=\frac{(n-s)_t}{(s+1)_t} \binom{n}{s}_{t}.\]
\subsection{Proof of Theorem~\ref{thm:Bsimply}}\label{subsec:fTh}
Firstly, we define iteratively the family of integers~$\{A^{k,n}_j\}$. Set $A^{k,n}_k=1$, for $h \in \{1, \dots, k-1\}$ we define 
\begin{equation}\label{def:A}
    A_h^{k, n} =	-\sum_{j=h+1}^{k}(-1)^{\lfloor\frac{j-h+1}{2}\rfloor} \binom{n-j+\lfloor\frac{j-h}{2}\rfloor}{\lfloor\frac{j-h}{2}\rfloor}A_{j}^{k,n}
\end{equation}
Moreover, by convention we set $ A_h^{k, n} =0$ if $h>k$ or if $h \leq 0$. Using properties of binomials and Equation~(\ref{def:A}) it is possible to prove that the integers $A^{k,n}_h$ satisfy 
nice iterative properties:
\begin{lemma}\label{lem:relA}
\begin{enumerate}
    \item $A_{h+1}^{k,n}=A_{h}^{k-1,n-1}$,
        \item $A_{h}^{k,k}=A_{h}^{k-1,k-1}+A_{h}^{k-2,k-1}$,
    \item $A_{h}^{k,n}=A_{h}^{k,n-1}+A_{h}^{k-2,n-1}$, if $k<n$.
    \end{enumerate}
\end{lemma}
We consider now the expression $\sum_{j=0}^{k}A_i^{k,n}R_i$.  It can be written in the form
\begin{equation}
 \Lambda^{k,n}_k C_k + \sum_{h=0}^{k-1} \Gamma_h^{k, n}C_h = 0,
\end{equation}
for some coefficients $\Gamma_h^{k, n}$ that we are going to determine explicitly. 
\begin{proposition}For every $h$ such that  $0 <h<k$ the equality
$\Gamma_h^{k, n} = \Gamma_0^{k-h, n-h}$ holds.
\end{proposition}
\proof
By definition we have that $\Gamma_h^{k,n}=\sum_{j=h}^k A_j^{k,n} \Lambda_h^{j,n}$. Now we use Lemma~\ref{lem:expansion} to expand~$\Lambda_h^{j,n}$:
\begin{align*}
    \Gamma_h^{k,n}&=\sum_{j=h+1}^k A_j^{k,n} \left[(-1)^{\lfloor\frac{j-h+1}{2}\rfloor} \binom{n-j+\lfloor\frac{j-h}{2}\rfloor}{\lfloor\frac{j-h}{2}\rfloor} \Lambda_h^{h,n}  + \Lambda_0^{j-h,n-h}\right] + A_h^{k,n} \Lambda_h^{h,n}\\
  &=\left[\sum_{j=h}^k (-1)^{\lfloor\frac{j-h+1}{2}\rfloor} \binom{n-j+\lfloor\frac{j-h}{2}\rfloor}{\lfloor\frac{j-h}{2}\rfloor}A_j^{k,n} \right] \Lambda_h^{h,n} + \sum_{j=h+1}^k A_j^{k,n} \Lambda_0^{j-h,n-h}\\
  \intertext{using Equation~(\ref{def:A}) and setting $t=j-h$ we have}
 &=\sum_{t=1}^{k-h} A_{t+h}^{k,n} \Lambda_0^{t,n-h} \\
 \intertext{and now by Lemma~\ref{lem:relA}}&= \sum_{t=1}^{k-h} A_{t}^{k-h,n-h} \Lambda_0^{t,n-h}=\Gamma_0^{k-h, n-h}.
\end{align*} \endproof
\begin{rmk}\label{rmk:G2}
By Equation~\ref{eqn:expansion2} we know that $\Lambda_0^{2,n}=-p(n,q,t)+\Lambda_0^{2,n-1}$. Moreover, observe that $\Lambda_0^{1,n}=\Lambda_0^{1,n-1}$ and $A_1^{2,n}=A_1^{2,n-1}=1$. We consequently obtain 
\[\Gamma_0^{2,n}=\Lambda_0^{2,n} + \Lambda_0^{1,n} = -p(n,q,t)+\Lambda_0^{2,n-1}+ \Lambda_0^{1,n-1}=  -p(n,q,t)+ \Gamma_0^{2,n-1}\]
\end{rmk}
\begin{proposition}\label{prop:expansion}If $k>2$, the following relations between the coefficients $\Gamma_0^{k,n}$ hold: 
\[\Gamma^{k,k}_0=\Gamma^{k-1,k-1}_0+\Gamma^{k-2,k-1}_0-\Gamma^{k-2,k-2}_0\] \[ \Gamma_0^{k,n}=\Gamma^{k,n-1}_0 - \Gamma^{k-2,n-2} + \Gamma^{k-2,n-1}_0 \qquad \mbox{ for } k<n\]
\end{proposition}
\proof
We consider $\Gamma_0^{k,n}=\sum_{j=1}^k A_j^{k,n} \Lambda_0^{j,n}$ end expand $\Lambda_0^{j,n}$ according to Lemma~\ref{lem:expansion}. We obtain 
\begin{align*}
    \Gamma_0^{k,n}&=\sum_{j=2}^k A_j^{k,n} \left[(-1)^{\lfloor\frac{j+1}{2}\rfloor} \binom{n-\lfloor\frac{j+1}{2}\rfloor -1}{\lfloor\frac{j}{2}\rfloor-1} p(n,q,t) - \Lambda_0^{j-2,n-2} + \Lambda_0^{j,n-1}\right] + A_1^{k,n} \Lambda_0^{1,n} \\
    &=\left[ \sum_{j=2}^k(-1)^{\lfloor\frac{j+1}{2}\rfloor} \binom{n-\lfloor\frac{j+1}{2}\rfloor -1}{\lfloor\frac{j}{2}\rfloor-1}A_j^{k,n} \right] p(n,q,t) - \sum_{j=3}^k A_j^{k,n}\Lambda_0^{j-2,n-2} + \sum_{j=1}^k A_j^{k,n}\Lambda_0^{j,n-1}\end{align*}
Observe now that Equation~(\ref{def:A}) implies 
 $$\sum_{j=2}^k(-1)^{\lfloor\frac{j+1}{2}\rfloor} \binom{n-\lfloor\frac{j+1}{2}\rfloor -1}{\lfloor\frac{j}{2}\rfloor-1}A_j^{k,n}=0.$$ Furthermore using Lemma~\ref{lem:expansion} and setting $t=j-2$ we have
    \begin{align*}
    \Gamma_0^{k,n}&=\sum_{j=1}^k A_j^{k,n}\Lambda_0^{j,n-1} - \sum_{j=3}^k A_j^{k,n}\Lambda_0^{j-2,n-2}\\
    &=  \sum_{j=1}^k \left[A_j^{k,n-1} + A_j^{k-2,n-1} \right]\Lambda_0^{j,n-1} - \sum_{t=1}^{k-2} A_{t+2}^{k,n}\Lambda_0^{t,n-2} \\
    &= \sum_{j=1}^k A_j^{k,n-1}\Lambda_0^{j,n-1} +  \sum_{j=1}^{k-2} A_j^{k-2,n-1}\Lambda_0^{j,n-1} - \sum_{t=1}^{k-2} A_{t}^{k-2,n-2}\Lambda_0^{t,n-2}\\
    &= \Gamma_0^{k,n-1}+\Gamma_0^{k-2,n-1} - \Gamma_0^{k-2,n-2}.
\end{align*}
and analogously 
\begin{align*}
    \Gamma_0^{k,k}&=\sum_{j=2}^k A_j^{k,k} \left[(-1)^{\lfloor\frac{j+1}{2}\rfloor} \binom{k-\lfloor\frac{j+1}{2}\rfloor -1}{\lfloor\frac{j}{2}\rfloor-1} p(k,q,t) - \Lambda_0^{j-2,k-2} + \Lambda_0^{j,k-1}\right] + A_1^{k,k} \Lambda_0^{1,k} \\
    &=\left[ \sum_{j=2}^k(-1)^{\lfloor\frac{j+1}{2}\rfloor} \binom{k-\lfloor\frac{j+1}{2}\rfloor -1}{\lfloor\frac{j}{2}\rfloor-1}A_j^{k,k} \right] p(k,q,t) - \sum_{j=3}^k A_j^{k,k}\Lambda_0^{j-2,k-2} + \sum_{j=1}^{k-1} A_j^{k,k}\Lambda_0^{j,k-1}\\
    &=  \sum_{j=1}^{k-1} \left[A_j^{k-1,k-1} + A_j^{k-2,k-1} \right]\Lambda_0^{j,k-1} - \sum_{t=1}^{k-2} A_{t}^{k-2,k-2}\Lambda_0^{j,k-2} \\
    &= \Gamma_0^{k-1,k-1}+\Gamma_0^{k-2,k-1} - \Gamma_0^{k-2,k-2}.
\end{align*}
\endproof
Making explicit computations for $n=2,3$, it is possible to prove that $\Gamma_0^{2,2}=\Gamma_0^{3,3}$. Moreover observe that  $\Gamma_0^{1,n}=\Gamma_0^{1,n+1}$ for every $n>1$. As a consequence of Proposition~\ref{prop:expansion} we obtain:
\begin{corollary}\label{cor:relations}
The following relations hold:
\[\Gamma^{k,k}_0= \begin{cases} \Gamma^{k-2,k-1}_0& \mbox{ if $k$ is even,}\\ \Gamma^{k-1,k-1}_0 & \mbox{ if $k$ is odd.}\end{cases}
\qquad 
\Gamma^{k,n}_0=
\begin{cases}
\Gamma^{k,n-1}_0 & \mbox{ if $n>k>2$ and $k$ is odd,} \\
\Gamma^{k-2,n-1}_0 & \mbox{ if $n>k>2$ and $k$ is even.}
\end{cases}
\]
\end{corollary}
Theorem~\ref{thm:Bsimply} comes directly by Remark~\ref{rmk:G2} and iterating the relations of Corollary~\ref{cor:relations}.
\subsection{Proof of Formulae in Type D}
We denote by $C_h(q,t)$ the formal series $C_{\omega_{2h}}(q,t)$. Moreover, if $n=2k+1$ (resp. $n=2k$) we denote by $C_k(q,t)$ the formal series $C_{\omega_{n-1}+\omega_n}(q,t)$ (resp. $C_{2\omega_{n}}(q,t)$). Our formulae can be obtained dealing with the non specialized version of Equation~4.4 of~\cite{SDT2}. The reduced recurrence $R_k$ for $C_k(q,t)$ can be written in the form 
\begin{equation}
R_k=\sum_{h \leq k}\Lambda_h^{k,n}(q,t)C_h(q,t) = 0.
\end{equation}
for certain coefficients $\Lambda_h^{k,n}(q,t)$.
Reasoning as in Remark~\ref{rem:Lambdak}, a non specialized analogue of Formula~4.5 in~\cite{SDT2} can be achieved: 
\begin{equation*}\label{coefficientek}
\Lambda_k^{k, \,n}(q,t)=\begin{cases}\frac{2(t^{2k}-1)(1-qt^{2k-1})}{t^{2k-1}(t-1)}& \mbox{if }    n=2k, \\\frac{(t^{2k}-1)(1-qt^{2(n-k)-1})}{t^{n-1}(t-1)}&\mbox{otherwise. }.\end{cases}
\end{equation*}
Set now
 \begin{equation*}
 \displaystyle{
b_{k,n} =
\left\{
	\begin{array}{lll}
\frac{(t-q)(t^{2k}-1)}{t^{k}(t-1)}& \mbox{if }    n=2k, \\
 	& \\
  \frac{(t-q)(t^{n}-1)(t^{n-2k}+1)}{t^{n-k}(t-1)}& \mbox{otherwise. } \\
	\end{array}
\right.}
\end{equation*}
The next Proposition is a non specialized version of Proposition 4.6 of \cite{SDT2}.
\begin{proposition}\label{prop:Dbella}The following recursive relation hold:
  \begin{equation}\label{Riduzionebella} \Lambda_k^{k, n}(q,t) {C}_{k}(q,t)- \sum_{i=1}^k b_{i,n-2(k-i)}{C}_{k-i}(q,t)=0\end{equation}
  \end{proposition}
  Using Equation~(\ref{Riduzionebella}) and Remark~\ref{rmk:base} as base step, it is possible to prove inductively Theorem~{\ref{thm:geD}}.
  \proof(of Theorem~{\ref{thm:geD}}) As observed in Remark~\ref{rem:ratioexp}, the formal series $\overline{C}_{k}(q,t)=C_k(q,t)/C_0(q,t)$ satisfies the Recurrence~(\ref{Riduzionebella}).
Specializing Equation~(\ref{Riduzionebella}) at $q=0$ and recalling that $E_{\omega_{2k}}(t)=\overline{C}_{k}(0,t)$, Formulae~(\ref{eqn:geD1}), (\ref{eqn:geD2}) and (\ref{eqn:geD3}) can be obtained recursively proving that
\begin{align*}\frac{(t^{2k}-1)}{t^{n-1}(t-1)} \frac{t^k(t^{n-2k}+1)}{(t^n+1)} \binom{n}{k}_{t^2}&= \frac{(t^{2k}-1)}{t^{n-1}(t-1)} \overline{C}_{k}(0,t)\\&= \sum_{i=1}^{k-1} \frac{t^{k+i}(t^{n-2i}-1)(t^{n-2k}+1)}{t^{n-1}(t-1)} \overline{C}_{i}(0,t)\\ &=\sum_{i=1}^{k-1} \frac{t^{k+i}(t^{n-2i}-1)(t^{n-2k}+1)}{t^{n-1}(t-1)} \frac{t^i(t^{n-2i}+1)}{(t^n+1)} \binom{n}{i}_{t^2}\\
&=\frac{t^k(t^{n-2k}+1)}{t^{n-1}(t-1)(t^n+1)} \sum_{i=1}^{k-1} t^{2i}(t^{2(n-2i)}-1) \binom{n}{i}_{t^2}.
\end{align*}
Again we reduced to Identity~(\ref{eqn:binomials}) that we just proved in Section~\ref{subsec:ReederB}. \endproof
\subsection{Proof of Proposition~\ref{prop:Dbella}} The proof is analogue to the proof of Proposition 4.6 contained in Section~5 of \cite{SDT2}. Set 
  \[r(n,q,t)=t^{(n-1)}-qt^{-(n-1)}+t^{-(n-2)}-qt^{(n-2)}=\frac{(t-q)\left(t^{2n-3}+1\right)}{t^{n-1}}\]
We denote by $\Omega_{0 }^{\lambda, \; n}$ the set of weights of the form $w\lambda$ conjugated to 0. If $\lambda=\omega_{2k}$ we will use the notation~$\Omega_{0 }^{k, \; n}$ Coherently with our previous notations, if $n=2k+1$ (resp. $n=2k$) we denote by $\Omega_{0 }^{k, \; n}$ the set ow weights of the form $w (\omega_{n-1}+\omega_n)$ (resp. $w(2\omega_n)$) conjugated to 0. We recall the following results by \cite{SDT2}, Section~4.1:
\begin{rmk}[\cite{SDT2}, Remark~4.3]\label{omegah+zero}
The weights giving non zero contribution to $\Lambda_{h}^{k,n}, k>h>0$ are of the form $e_1+\dots+e_{2h}+\nu$, where $\nu$ has the first $2h$ coordinates equal to 0. Considering the immersion of $D_{n-2h} \rightarrow D_n$ induced by the Dynkin diagrams, $\nu$ can be then identified with a weight in $\Omega_{0}^{k-h,n-2h}$.
\end{rmk}
\begin{lemma}[\cite{SDT2}, Lemma~4.5]\label{zerocontribution2}
Set $\lambda=\omega_{2k}$, $2k<n$ or $\lambda=\omega_{n-1}+\omega_{n}$, $n=2k+1$ and let $w\in W$ be such that $w \lambda$ is conjugated to $0$. Then $w \lambda$ is of one of the following form:
\begin{enumerate}
 \item The $2k$ non zero coordinates of $w \lambda$ are pair of consecutive coordinates $((w\lambda)_{(j)},(w\lambda)_{(j)+1})$  of the form $(-1,1)$.
 \item There are $2(k-1)$ non zero coordinates that are pair of consecutive coordinates $((w\lambda)_{(j)},(w\lambda)_{(j)+1})$  of the form $(-1,1)$ and the latter two are equal to $-1$.
\end{enumerate}
In both cases there exists an element $\sigma \in W$ of length $l(\sigma)= k$ such that $\sigma (w \lambda + \rho)= \rho$.
\end{lemma}
\begin{rmk}
Consider $\mu \in \Omega_{0}^{k, \, n}$ and denote by $\epsilon_n$ the sign change on the $n$-th coordinate.
\begin{itemize}
\item \emph{If $n=2k$} then $\mu$ must be of the form $-e_1+e_2+\nu$ with $\nu \in \Omega_0^{k-1,2k-2}$,
\item \emph{If $n=2k+1$} then $\mu$ must be of the form $-e_1+e_2+\nu$ with $\nu \in \Omega_0^{k-1,2k-1}$ or $\mu=(0, \mu_2, \dots, \mu_n)$ where $\mu'=(\mu_2, \dots, \mu_n) \in \Omega_0^{k,2k}$ or $\epsilon_n \mu' \in \Omega_0^{k,2k}$,
\item \emph{ If $n\neq 2k, 2k+1$} then $\mu$ must be of the form $-e_1+e_2+\nu$ with $\nu \in \Omega_0^{k-1,n-2}$ or $\mu=(0, \mu_2, \dots, \mu_n)$ where $\mu'=(\mu_2, \dots, \mu_n) \in \Omega_0^{k,n-1}$. 
\end{itemize}
\end{rmk}
The above considerations lead to non specialized versions of recursive relations (4.6), (4.7), (4.8) and (4.9) in \cite{SDT2}:
\begin{equation}\label{Lh}
\Lambda_h^{k, \, n}(q,t)= (-1)^{k-h}\Lambda_h^{h,n}(q,t) |\Omega_0^{k-h, n-2h}| + \Lambda_0^{k-h, \, n-2h}(q,t).  \end{equation}
\begin{equation}\label{L002k}
 \Lambda_0^{k,2k}(q,t) = (-1)^{k}\sum_{i=1}^{k} r(2i, q,t)=(-1)^{k} r(2k,q,t) - \Lambda_0^{k-1, 2k-2}(q,t).
\end{equation}
 \begin{equation}\label{L02k1}\Lambda_{0}^{k, \, 2k+1}(q,t)= (-1)^k r(2k+1,q)|\Omega_0^{k-1,2k-1}|- \Lambda_0^{k-1,2k-1}(q,t) + 2 \Lambda_0^{k,2k}(q,t), \end{equation}
  \begin{equation}\label{L0n}\Lambda_{0}^{k, \, n}(q,t)= (-1)^k r(n,q,t)|\Omega_0^{k-1,n-2}|- \Lambda_0^{k-1,n-2}(q,t) + \Lambda_0^{k,n-1}(q,t). \end{equation}
  where 
  \begin{equation*}\label{CardG0}
 \displaystyle{
| \Omega_0^{k, n}| =
\left\{
	\begin{array}{lll}
	1 & \mbox{if } \lambda = 2 \omega_n \mbox{ or if }\lambda=0 \\
    \frac{n}{k}\binom{n-k-1}{k-1}& \mbox{if } \lambda=\omega_{2k}  \mbox{ and }   2k < n \mbox{ or } \lambda=\omega_{n-1}+\omega_{n}  \mbox{ and } n=2k+1.\\
	\end{array}
\right.}
\end{equation*}
As in Section~5 of \cite{SDT2}, define a family of integers $A_h^{k.n}$
in the following way:
\begin{equation}\label{blabla}
A_h^{k, n} =
\left\{
	\begin{array}{lll}
        0  & \mbox{if } h>k \mbox{ or } h\leq 0,\\
		1  & \mbox{if } h=k, \\
		-\sum_{i=h+1}^k(-1)^{i-h} | \Omega_0^{i-h, n-2h}| A_{i}^{k,n} & \mbox{otherwise.}
	\end{array}
\right.
\end{equation}
and consider 
 \begin{equation}\label{RiduzionebellaS5}\sum_{i=1}^k A_i^{k,n}R_i = \Lambda_k^{k, n}(q,t) C_{k}(q,t)- \sum_{i=0}^{k-1} \Gamma_i^{k,n}(q,t)C_{i}(q,t).\end{equation}
Performing the same computation as in Proposition~5.2 of \cite{SDT2} it is possible to prove that $\Gamma_h^{k,n}(q,t)=\Gamma_0^{k-h,n-2h}(q,t)$ if $k>h>0$. Analogously the following relations hold: 
 \[\Gamma_0^{k,2k}(q,t)=-\sum_{j=2}^{k+1}r(j,q,t), \qquad \Gamma_0^{k,2k+1}(q,t)=2 \Gamma_0^{k,2k}(q,t)-r(k+2,q,t), \]\[ \Gamma_0^{k,n}(q,t)= \Gamma_0^{k,n-1}(q,t)-r(n-k+1,q,t).\]
 Now it is straightforward to show that $\Gamma_0^{k,n}(q,t)=b_{k,n}$ and 
 then~$\Gamma_0^{k-h,n-2k}(q,t)=b_{k-h,n-2h}$.

\subsection{Open Questions about Generalized Exponents and Small Representations.}\label{subsec:conj}
Some natural questions arise as consequences of our results. Firstly, as a consequence of Theorem~\ref{Broer}, generalized exponents of small representations are related to the so called \emph{fake degrees}, i.e. degrees of generators (as $S(\h)^W$-module) of isotypic components of $W$-representations in~$S(\h)$. There exists an ample literature about fake degrees and many formulae to obtain them in terms of suitable combinatorial statistics (see \cite{Opd} for a complete survey about the topic and \cite{ABR} \cite{BC}, \cite{BC1}, \cite{GP}, \cite{KW}, \cite{RY}, \cite{St2} for more specific results). It could be interesting to find a purely combinatorial proof of Formulae (\ref{eqn:geC}), (\ref{eqn:geB}), (\ref{eqn:geD1}), (\ref{eqn:geD2}) and (\ref{eqn:geD3}).\\
Moreover, a closer analysis of formulae proved in \cite{SDT} and \cite{SDT2} for graded multiplicities of small representations in the exterior algebra, underlines some similarities with the results contained in \cite{DCMPP} and \cite{DCPP}.
In fact, Theorem \ref{thm:DCPP} and its analogous version for little adjoint representation are suggested by the following factorizations of $P(\g, \Lambda \g, q)$ and $P(V_{\theta_s}, \Lambda \g, q)$:
  \begin{equation}\label{eqn:facagg}P(\g, \Lambda \g, q)= (1+q^{-1})\prod_{i=1}^{n-1}(q^{2e_i+1}+1) {E_\theta( q^2)}  \end{equation}
     \begin{equation}\label{eqn:faclagg}P(V_{\theta_s}, \Lambda \g, q)= (1+q^{-1})\prod_{i=1}^{n-1}(q^{2e_i+1}+1) {E_{\theta_s}( q^2)}  \end{equation}
In particular, the authors of \cite{DCMPP} and \cite{DCPP} noticed that the factor $\prod_{i=1}^{n-1}(q^{2e_i+1}+1)$ is the Poincarè polynomial of the exterior algebra over the first $n-1$ generators $P_1, \dots , P_{n-1}$ of the algebra of invariants in $\Lambda \g$.
A direct computation shows that similar factorizations can be achieved for polynomials of graded multiplicities of certain small representations. As an example, comparing the results exposed in Theorem~\ref{thm:geB} with formulae proved in \cite{SDT}, in type~$B_n$ the polynomials for graded multiplicities can be rearranged as
\[
 P(V_{\omega_{2s}}, \Lambda \g, q)=(1+q^{-1})\prod_{i=1}^{n-s}(1+q^{2e_i+1})
\prod_{i=1}^{s-1}(1+q^{2e_i+1})E_{\omega_{2s}}( q^2),
\]
\[
 P(V_{\omega_{2s+1}}, \Lambda \g, q)=(1+q^{-1})\prod_{i=1}^{s}(1+q^{2e_i+1})
\prod_{i=1}^{n-s-1}(1+q^{2e_i+1})E_{\omega_{2s+1}}( q^2).
\] Analogously, using Theorem~\ref{thm:geC}, in type $C_n$ it is possible to obtain the factorization:
 \[P(V_{\omega_{2k}}, \Lambda \g,  q)= (1+q^{-1})\prod_{i=1}^{n-k}(q^{2e_i+1}+1)\prod_{i=1}^{k-1}(q^{2e_i+1}+1)E_{\omega_{2k}} (q^2). \]
Consequently, it is natural to ask if there exist examples of small representations $V_\lambda$, different from $V_\theta$ and $V_{\theta_s}$, such that the module $\Hom_\g(V_\lambda, \Lambda \g)$ has a structure of free module over a suitable exterior algebra of invariants,
with degrees prescribed by factorizations of  $P(V_{\lambda}, \Lambda \g, q)$ similar to the ones in Formulae (\ref{eqn:facagg}) and (\ref{eqn:faclagg}). 

\end{document}